\documentclass[12pt]{article}

\usepackage{indentfirst}
\usepackage{amsfonts}
\usepackage{amsmath}
\usepackage{amssymb}
\usepackage{amsbsy, amsthm}
\usepackage[margin=1in]{geometry}

\newtheorem{dfn}{Definition}[section]
\newtheorem{theorem}[dfn]{Theorem}
\newtheorem{lemma}[dfn]{Lemma}
\newtheorem{proposition}[dfn]{Proposition}

\newtheorem{corollary}[dfn]{Corollary}
\newtheorem{conjecture}[dfn]{Conjecture}
\newtheorem{question}[dfn]{Question}
\newenvironment{pf}{\noindent{\bf Proof.}}
{\enspace\vrule height5pt depth0pt width5pt}

\def\X {{\mathcal X}}

\def\Q {{\mathcal Q}}

\def\W {{\mathcal W}}
\def\adh {{\rm adh}}
\def\C {{\mathcal C}}

\begin{document}
\title{Immersion and clustered coloring}
\author{Chun-Hung Liu\thanks{chliu@math.tamu.edu. Partially supported by NSF under Grant No.~DMS-1929851 and DMS-1954054.} \\
\small Department of Mathematics, \\
\small Texas A\&M University,\\
\small College Station, TX 77843-3368, USA}

\maketitle

\begin{abstract}
\noindent
Hadwiger and Haj\'{o}s conjectured that for every positive integer $t$, $K_{t+1}$-minor free graphs and $K_{t+1}$-topological minor free graphs are properly $t$-colorable, respectively.
Clustered coloring version of these two conjectures which only require monochromatic components to have bounded size has been extensively studied. 
In this paper we consider the clustered coloring version of the immersion-variant of Hadwiger's and Haj\'{o}s' conjecture proposed by Lescure and Meyniel and independently by Abu-Khzam and Langston.
We determine the minimum number of required colors for $H$-immersion free graphs, for any fixed graph $H$, up to a small additive absolute constant.
Our result is tight for infinitely many graphs $H$.

A key machinery developed in this paper is a lemma that reduces a clustering coloring problem on graphs to the one on the torsos of their tree-cut decomposition or tree-decomposition.
A byproduct of this machinery is a unified proof of a result of Alon, Ding, Oporowski and Vertigan and a result of the author and Oum about clustered coloring graphs of bounded maximum degree in minor-closed families.
\end{abstract}

\section{Introduction}
All {\it graphs} in this paper are finite and allowed to have loops and parallel edges.
Graph coloring is one central research direction in graph theory.
For a positive integer $t$, a graph is {\it properly $t$-colorable} if it can be (vertex-)partitioned into $t$ edgeless induced subgraphs. 
Every properly $t$-colorable graph does not contain $K_{t+1}$ as a subgraph.
But the converse statement is not true: for every graph $H$ that contains a cycle, there exists no integer $C$ such that every graph with no $H$ subgraph is properly $C$-colorable \cite{e}.

A line of research focuses on coloring graphs that forbid a fixed graph $H$ as a more general structure than subgraphs.
One of the most famous problems in this direction is a conjecture of Hadwiger \cite{h} stating that for every positive integer $t$, every graph with no $K_{t+1}$-minor\footnote{A graph $G$ contains an {\it $H$-minor} for some graph $H$ if $H$ is isomorphic to a graph that can be obtained from a subgraph of $G$ by contracting edges.} is properly $t$-colorable.
Hadwiger's conjecture is very difficult, and the case $t=4$ implies the celebrated Four Color Theorem \cite{ah,ahk,rsst}.
Hadwiger's conjecture is true for $t \leq 5$ \cite{h,rst,w_color} and open for $t \geq 6$.
Norin and Song \cite{ns} recently improved an old general upper bound given independently by Kostochka \cite{k1,k2} and Thomason \cite{t}; more recently, Postle \cite{ns} joined the project to further improve the result. 
Even more recently, Delcourt and Postle \cite{dp} further improved the result by proving that every graph with no $K_t$-minor is properly $O(t\log\log t)$-colorable.

A similar conjecture was proposed by Haj\'{o}s in 1940s stating that for every positive integer $t$, every graph with no $K_{t+1}$-topological minor\footnote{A graph $G$ contains an {\it $H$-topological minor} for some graph $H$ if some subgraph of $G$ is isomorphic to a subdivision of $H$.} is properly $t$-colorable.
Haj\'{o}s' conjecture is stronger than Hadwiger's conjecture and is true for $t \leq 3$ \cite{d_hajos}.
But it is too strong to be true in general: Catlin \cite{c_hajos} disproved the cases for $t \geq 6$, and Erd\H{o}s and Fajtlowicz \cite{ef} proved that $\Omega(t^2/\log t)$ colors are required.
The cases for $t \in \{4,5\}$ remain open.

Due to the difficulty of Hadwiger's conjecture and the incorrectness of Haj\'{o}s' conjecture, relaxations of those two conjectures have been extensively considered.
One relaxation is to consider clustered coloring.

For positive integers $t$ and $N$, we say that a graph $G$ is {\it $t$-colorable with clustering $N$} if $G$ can be (vertex)-partitioned into $t$ induced subgraphs with no component on more than $N$ vertices. 

For every positive integer $t$, define $f(t)$ (and $g(t)$, respectively) to be the minimum $k$ such that there exists an integer $N$ such that every graph with no $K_{t+1}$-minor (and $K_{t+1}$-topological minor, respectively) is $k$-colorable with clustering $N$.
Clearly, $g(t) \geq f(t)$.
Edwards, Kang, Kim, Oum and Seymour \cite{ekkos} showed that $f(t) \geq t$ by using essentially the same method of Linial, Matou\v{s}ek, Sheffet and Tardos \cite{lmst} who proved $f(t) \geq t-1$, so the number of required colors for the clustered coloring version of Hadwiger's conjecture and Haj\'{o}s' conjecture are the same as their original version.
The clustered coloring version of Hadwiger's conjecture has been extensively studied \cite{dn,ekkos,hw,km,lo,lw_minor,lw_topo,w_clu}.
In particular, the author and Wood \cite{lw_topo} proved that for every integer $t$ and every graph $H$, there exists an integer $N$ such that every graph with no $K_{t+1}$-topological minor and $H$-minor is $(t+1)$-colorable with clustering $N$.
This implies $f(t) \leq t+1$ which is the currently best known upper bound in the literature\footnote{Dvo\v{r}\'{a}k and Norin \cite{dn} announced that a forthcoming paper will prove $f(t)=t$. This result will solve the clustered coloring version of Hadwiger's conjecture. But it is incomparable with the aforementioned result in \cite{lw_topo}.}.
For the clustered coloring version of Haj\'{o}s' conjecture, the author and Wood \cite{lw_topo} proved that $g(t) \leq \max\{4t-5,1\}$ which is the only known linear upper bound in the literature.

Another well-known coloring problem about forbidding a complete graph as a more general structure is stated in terms of immersions.
For two distinct edges $e_1,e_2$ of a graph $G$ with a common end $v$, {\it splitting off $e_1$ and $e_2$ along $v$} is the operation that deletes $e_1$ and $e_2$ and adds an edge $(e_1-\{v\}) \cup (e_2-\{v\})$.
For a graph $H$, we say that a graph $G$ {\it contains an $H$-immersion} if $H$ is isomorphic to a graph that can be obtained from a subgraph of $G$ by repeatedly splitting off edges and deleting isolated vertices.
Clearly, if $G$ contains an $H$-topological minor, then $G$ contains an $H$-immersion and an $H$-minor.
But the minor relation is incomparable with the immersion relation.
Lescure and Meyniel \cite{lm} and Abu-Khzam and Langston \cite{al} independently proposed an immersion version of Hadwiger's conjecture. 

\begin{conjecture}[\cite{al,lm}] \label{imm_conj}
For every positive integer $t$, every graph with no $K_{t+1}$-immersion is properly $t$-colorable.
\end{conjecture}

The cases for $t \leq 3$ of Conjecture \ref{imm_conj} follow from the correctness of Haj\'{o}s' conjecture for $t \leq 3$.
DeVos, Kawarabayashi, Mohar and Okamura \cite{dkmo} proved the cases $4 \leq t \leq 6$.
It is open for $t \geq 7$.
The general upper bound for Conjecture \ref{imm_conj} has been steadily improved \cite{ddfmms,dy,glw}, and the currently best upper bound is $3.54t+7.54$ due to Gauthier, Le and Wollan \cite{glw}. 

Unlike the clustered coloring version of Hadwiger's and Hajos' conjectures, it was unknown whether it requires $t$ colors to color graphs with no $K_{t+1}$-immersion with bounded clustering.
The main result of this paper shows that $t$ colors are indeed much more than needed.
Our main result is actually more general and is about graphs with no $H$-immersion for any fixed graph $H$.

Let $H$ be a graph with maximum degree $d$.
As graphs with maximum degree at most $d-1$ cannot contain an $H$-immersion, the number of required colors for graphs with maximum degree at most $d-1$ gives a lower bound for the number of required colors for graphs with no $H$-immersion.
Our main result shows that this lower bound is very close to the correct value.

Define $\chi_*: {\mathbb N} \cup \{0\} \rightarrow {\mathbb N}$ to be the function such that for every $x \in {\mathbb N} \cup \{0\}$, $\chi_*(x)$ is the minimum $k$ such that there exists $N_x \in {\mathbb N}$ such that every graph of maximum degree at most $x$ is $k$-colorable with clustering $N_x$.
Note that $\chi_*$ exists and $\chi^*(x) \leq x+1$ for every $x \in {\mathbb N} \cup \{0\}$, since every graph of maximum degree at most $x$ is properly $(x+1)$-colorable.
Haxell, Szab\'{o} and Tardos \cite{hst} proved that $\chi_*(x) \leq \lceil \frac{x+1}{3} \rceil$ for every $x \in {\mathbb N} \cup \{0\}$, and there exist $\epsilon>0$ and $D$ such that $\chi_*(x) \leq (\frac{1}{3}-\epsilon)x$ for every $x \geq D$.
On the other hand, it is known that $\chi_*(x) \geq \lfloor \frac{x+6}{4} \rfloor$ for every $x \in {\mathbb N}$ \cite{adov,hst}.
However, even the asymptotic behavior of $\frac{\chi_*(x)}{x}$ remains unknown.

For a graph $H$, define $\chi_*(H)$ to be the minimum $k$ such that there exists $N \in {\mathbb N}$ such that every graph with no $H$-immersion is $k$-colorable with clustering $N$.
The following is the main theorem of this paper.

\begin{theorem} \label{main_intro}
Let $d$ be a positive integer, and let $H$ be a graph of maximum degree $d$.
	\begin{enumerate}
		\item If $d=1$, then $\chi_*(H)=1$.
		\item If $d \geq 2$ and $H$ has exactly one vertex of degree $d$, then $\chi_*(d-1) \leq \chi_*(H) \leq \max\{\chi_*(d-2)+1,4\}$.
		\item If $d \geq 2$ and $H$ has at least two vertices of degree $d$, then $\chi_*(d-2)+1 \leq \chi_*(H) \leq \max\{\chi_*(d-1)+1,4\}$.
	\end{enumerate}
\end{theorem}

Note that for every $x \in {\mathbb N}$, every graph with maximum degree at most $x$ can be partitioned into a stable set and an induced subgraph with maximum degree at most $x-1$, so $\chi_*(x) \leq \chi_*(x-1)+1$.
Since $\lfloor \frac{x+6}{4} \rfloor \leq \chi_*(x) \leq \lceil \frac{x+1}{3} \rceil$ for every positive integer $x$, there are infinitely many positive integers $d$ such that $\chi_*(d-1)=\chi_*(d-2)+1$, and there are infinitely many positive integers $d$ such that $\chi_*(d-2) = \chi_*(d-1)$.
Hence each statement of Theorem \ref{main_intro} is tight for infinitely many positive integers $d$.

We remark that some readers might get an impression at first glance that Theorem \ref{main_intro} is not surprising because other work of the author \cite{l_imm_str} has shown that $H$-immersion free graphs can be ``decomposed'' into graphs that are more or less close to graphs of maximum degree less than the maximum degree of $H$.
(See Theorem \ref{global} for a precise description.)
However, this impression is not correct in general. 
One example is that there exists no integer $k$ such that graphs of tree-width at most $w$ are $k$-colorable with bounded clustering for all integers $w$, even though graphs of tree-width at most $w$ can be ``decomposed'' in a similar flavor into graphs on at most $w+1$ vertices, which are 1-colorable with bounded clustering.
Other examples include the known results about finding nearly tight bounds for clustered chromatic number.
Those results\footnote{Such as in \cite{lw_minor} or the work of Dvo\v{r}\'{a}k and Norin about clustered coloring version of Hadwiger's conjecture.} are on graphs that are known to be ``decomposable'' into graphs whose clustered chromatic number are known, but they are still very difficult.

One contribution of this paper is to show that the decomposition used in \cite{l_imm_str} is nicely compatible with clustered coloring so that the generally false impression for decomposition is correct for this setting.
The following theorem is a key lemma for proving Theorem \ref{main_intro} and shows how to construct a clustered coloring for a graph with a given tree-cut decomposition of bounded bag size and bounded adhesion, where the number of colors is the required number of colors for coloring its torsos. 
(Notions related to tree-cut decomposition will be defined in Section \ref{sec:def}.)

\begin{theorem} \label{tree_cut_intro}
For any positive integers $N,\eta$ and $\alpha$, there exists a positive integer $N^*$ such that the following holds.
Let $G$ be a graph that admits a tree-cut decomposition $(T,\X=(X_t: t \in V(T)))$ of adhesion at most $\eta$ such that every bag contains at most $\alpha$ vertices.
For every $t \in V(T)$, let $k_t$ be a positive integer with $k_t+\lvert X_t \rvert \geq 2$ such that the graph obtained from the torso at $t$ by deleting $X_t$ is $k_t$-colorable with clustering $N$.
Then $G$ is $\max_{t \in V(T)}\{k_t + \min\{\lvert X_t \rvert,1\}\}$-colorable with clustering $N^*$. 
\end{theorem}

Theorem \ref{tree_cut_intro} leads to the following corollary showing that one can reduce clustered coloring problems on graphs of bounded maximum degree to the ones on the torsos of its tree-decomposition with bounded adhesion.
(Notions related to tree-decomposition will be defined in Section \ref{sec:tree_decomp}.)

\begin{corollary} \label{td_intro}
For any positive integers $\eta,d$ and $N$, there exists a positive integer $N^*$ such that the following holds.
For every integer $k$ and graph $G$ with maximum degree at most $d$ admitting a tree-decomposition $(T,\X)$ of adhesion at most $\eta$,
	\begin{enumerate}
		\item if $k\geq 2$ and for every $t \in V(T)$, the torso at $t$ in $(T,\X)$ is $k$-colorable with clustering $N$, then $G$ is $k$-colorable with clustering $N^*$, and
		\item if for every $t \in V(T)$, the subgraph of $G$ induced by the bag at $t$ is $k$-colorable with clustering $N$, then $G$ is $(k+1)$-colorable with clustering $N^*$.
	\end{enumerate}
\end{corollary}

Statement 1 in Corollary \ref{td_intro} leads to a simple unified proof of a result of Alon, Ding, Oporowski and Vertigan \cite{adov} and a result of the author and Oum \cite{lo}. 
The former states that graphs of bounded maximum degree and bounded tree-width are 2-colorable with bounded clustering and immediately follows from Corollary \ref{td_intro}.
The latter states that for every graph $H$, $H$-minor free graphs with bounded maximum degree are 3-colorable of bounded clustering.
Such graphs have tree-decompositions of bounded adhesion such that each torso can be made a graph of bounded layered tree-width by deleting a bounded number of vertices \cite{dmw}.
So the torsos are 3-colorable with bounded clustering \cite{lw_layer}, and hence so are the entire graphs by Statement 1 in Corollary \ref{td_intro}.

We remark that tree-cut decomposition and tree-decomposition of graphs are equivalent to expressions of graphs as edge-sums and clique-sums of their torsos, respectively.
So Theorem \ref{tree_cut_intro} and Corollary \ref{td_intro} reduce clustered coloring problems on graphs to the ones on the summands of edge-sums and clique-sums and hence are expected to have further applications.
Note that the bounded maximum degree condition for the clique-sum case is required, as there exists no $k$ such that for every $w$, graphs of tree-width at most $w$ is $k$-colorable with bounded clustering.

This paper is organized as follows.
In Section \ref{sec:def}, we include some necessary definitions.
We prove Theorem \ref{tree_cut_intro} in Section \ref{sec:tree-cut}.
In Section \ref{sec:proof}, we show how to combine Theorem \ref{tree_cut_intro} and work in \cite{l_imm_str} to prove Theorem \ref{main_intro}. 
We deduce Corollary \ref{td_intro} from Theorem \ref{tree_cut_intro} in Section \ref{sec:tree_decomp}.
We include some concluding remarks in Section \ref{sec:remarks}.

\section{Notations} \label{sec:def}

Let $G$ be a graph.
A {\it tree-cut decomposition} of $G$ is a pair $(T,\X)$ such that $T$ is a tree, and $\X$ is a collection $(X_t: t \in V(T))$ of pairwise disjoint (not necessarily non-empty) subsets of $V(G)$ such that $\bigcup_{t \in V(T)}X_t = V(G)$.
In addition,
	\begin{itemize}
		\item for every $t \in V(T)$, the set $X_t$ is called the {\it bag} at $t$;
		\item for every subset $S$ of $V(T)$, we define $X_S$ to be $\bigcup_{t \in S}X_t$; for every subgraph $S$ of $T$, we define $X_S$ to be $\bigcup_{t \in V(S)}X_t$;
		\item for any edge $t_1t_2$ of $T$, the {\it adhesion set} of $t_1t_2$ in $(T,\X)$, denoted by $\adh_{(T,\X)}(t_1t_2)$, is the set of edges of $G$ with one end in $X_{T_1}$ and one end in $X_{T_2}$, where $T_1$ and $T_2$ are the components of $T-t_1t_2$;
		\item the {\it adhesion} of $(T,\X)$ is $\max_{e \in E(T)} \lvert \adh_{(T,\X)}(e) \rvert$;
		\item the {\it torso} at $t$ in $(T,\X)$ is the graph obtained from $G$ by, for each edge $e$ of $T$ incident with $t$, identifying $X_{T_{t,t'}}$ into a vertex and deleting all loops incident with this new vertex, where $t'$ is the end of $e$ other than $t$, and $T_{t,t'}$ is the component of $T-t$ containing $t'$;
		\item each vertex in the torso at $t$ but not in $X_t$ is called a {\it peripheral vertex}.
	\end{itemize}
Note that for every $t \in V(T)$ and every edge $e$ of the torso at $t$, $e$ corresponds to an edge $e'$ of $G$ such that there exists no component $T'$ of $T-t$ such that $X_{T'}$ contains all ends of $e'$.

Let $G$ be a graph.
An {\it edge-cut} $[A,B]$ of a graph $G$ is a pair of disjoint subsets of $V(G)$ such that $A \cup B = V(G)$.
The {\it order} of an edge-cut $[A,B]$ of $G$ is the number of edges of $G$ with one end in $A$ and one end in $B$.

Let $(T,\X)$ be a tree-cut decomposition of a graph $G$.
Let $t$ be a node of $T$ or a connected subgraph of $T$.
Let $e$ be an edge of $T$ with at most one end in $t$.
We define $[A_{e,t},B_{e,t}]$ to be the edge-cut of $G$ with $B_{e,t} = \bigcup_{t''} X_{t''}$, where the union is over all nodes $t''$ contained in the component of $T-e$ containing $t$.

Let $G$ be a graph and let $S$ be a subset of $V(G)$.
We define $G[S]$ to be the subgraph of $G$ induced by $S$.
For each vertex $v$ of $G$, the {\it degree} of $v$ is the number of edges of $G$ incident with $v$, where each loop is counted twice.
The {\it maximum degree} of $G$ is the maximum of a degree of a vertex of $G$.

We say that a subgraph $H$ of a graph $G$ is {\it incident} with an edge $e$ of $G$ if $V(H)$ contains at least one end of $e$.

\section{From torsos to the whole graph} \label{sec:tree-cut}

The objective of this section is proving Theorem \ref{tree_cut_intro}.
The main challenge of the proof lies in the special case that every bag in the tree-cut decomposition has at most 1 vertex.
Theorem \ref{tree_cut_intro} follows from this special case easily, as shown in the proof of Lemma \ref{coloring_const_bag}.
This special case will be proved in Lemma \ref{coloring_bag1}, and we provide a sketch the proof of this special case before we formally prove it.

Let $(T,\X)$ be a tree-cut decomposition of a graph $G$ with $\lvert X_t \rvert \leq 1$ for every $t \in V(T)$, where $\X=(X_t: t \in V(T))$.
We find a depth-first-search ordering $t_1,t_2,...,t_{\lvert V(T) \rvert}$ of nodes of $T$.
Then starting from $i=1$ until $i=\lvert V(T) \rvert$, we greedily color every uncolored vertex that is either in $X_{t_i}$ or incident with an edge of $G$ corresponding to an edge of the torso at $t_i$ by a certain rule.
Note that at any time during the process, we can order the current monochromatic components by saying that a monochromatic component $C_1$ is older than another monochromatic component $C_2$ if the earliest colored vertex in $C_1$ is colored earlier than the earliest colored vertex in $C_2$.
In the formal proof, we will define the ordering of monochromatic components more carefully and describe it formally by using the notation $\sigma(C)$. 

Now we explain the rule for assigning colors at the time $i$ in the greedy algorithm.
When $X_{t_i}$ is nonempty and the vertex in $X_{t_i}$ is uncolored, we color this vertex with a color that is different from the color of the oldest monochromatic component adjacent to it to stop the growth of the oldest monochromatic component.
(If no such monochromatic component exists, then color the vertex arbitrarily.)
By the assumption, there exists a coloring $c_i'$ of the graph obtained from the torso at $t_i$ by deleting $X_{t_i}$ with bounded clustering by using $k_t$ colors in $[k_t+\lvert X_{t_i} \rvert]$ without using the color that was just used for coloring $X_{t_i}$.
Note that each uncolored vertex $v$ incident with an edge of $G$ corresponding to an edge of the torso at $t_i$ is contained in the union of bags of nodes in a component of $T-t_i$ corresponding to a peripheral vertex in the torso at $t_i$, and we call the color in $c_i'$ of this peripheral vertex the ``default color'' for $v$.
If $v$ is adjacent to some current monochromatic component, then we color $v$ by using a color different from the color of the oldest monochromatic component adjacent to $v$; otherwise, we color $v$ by using its default color.
This essentially completes the description of our greedy algorithm, except that we will also color some ``special vertices'' during the algorithm due to some technical reasons that will not be described in the proof sketch.

The main challenge is to bound the size of each final monochromatic component $C$.
The size of $C$ is determined by two factors: ``breadth'' and ``depth''.
The breadth of $C$ counts the number of components in $T-t_i$ whose bags intersect $C$, where $i$ is the first time in the algorithm such that some vertex in $C$ is colored.
The depth of $C$ counts the number of vertices in $C$ contained in the union of the bags contained in a single component of $T-t_i$.
The size of $C$ is bounded if both the breadth and the depth of $C$ are bounded.
We will prove that the breadth of $C$ is essentially bounded by the clustering of $c_i'$ and the size of the monochromatic components older than $C$ adjacent to $C$; and the depth is essentially bounded by the size of the monochromatic components older than $C$ adjacent to $C$.
We will also show that the number of those older components can be bounded by the adhesion of $(T,\X)$.
So the size of $C$ is bounded and the proof sketch is completed.

Claims 1-3 in the proof of Lemma \ref{coloring_bag1} are dedicated to showing the consistency of the ordering of monochromatic components during the algorithm.
Claims 4-7 are dedicated to bounding the depth of a monochromatic component.
Claims 8 and 9 are dedicated to bounding the breadth of a monochromatic component.

Now we introduce some terminologies that will be used in our proofs.
Let $k$ and $N$ be positive integers.
A {\it $k$-coloring} of a graph $G$ is a function $f:V(G) \rightarrow [k]$.
For a $k$-coloring $c$ of $G$, a {\it $c$-monochromatic component} (or a {\it monochromatic component in $c$}) is a component of $G[c^{-1}(\{i\})]$ for some $i \in [k]$.
For a graph $G$ and a function $f$ whose domain is a subset $S$ of $V(G)$, we say a vertex $v$ of $G$ is {\it $f$-colored} if $v \in S$, and we say $v$ is {\it $f$-uncolored} if $v \not \in S$.

\begin{lemma} \label{coloring_bag1}
For any positive integers $N$ and $\xi$, there exists a positive integer $N^*=N^*(N,\xi)$ such that the following holds.
Let $G$ be a graph that admits a tree-cut decomposition $(T,\X=(X_t: t\in V(T)))$ of adhesion at most $\xi$ such that every bag contains at most 1 vertex.
For every $t \in V(T)$, let $k_t$ be a positive integer such that $k_t +\lvert X_t \rvert \geq 2$ and the graph obtained from the torso at $t$ by deleting $X_t$ is $k_t$-colorable with clustering $N$.
Then $G$ is $\max_{t \in V(T)}\{k_t + \lvert X_t \rvert\}$-colorable with clustering $N^*$. 
\end{lemma}

\begin{pf}
Let $N$ and $\xi$ be positive integers.
Define the following.
\begin{itemize}
	\item Let $N_0 = \xi^2+\xi$.
	\item Let $N_1 = 2N_0\xi$.
	\item Let $N_2 = (1+2\xi^2(\xi+1)N)N_1^2 + N_0$.
	\item Let $f$ be the function with domain ${\mathbb N}$ such that 
		\begin{itemize}
			\item $f(1)=N_0$,
			\item for every $x \in {\mathbb N}$, $f(x+1)=(\xi+1) N_2 \cdot \sum_{i=1}^{x}f(i)$.
		\end{itemize}
	\item Define $N^*=1+(1+N\xi)(\xi+1)N_0 \cdot f(\xi)$.
\end{itemize}

Let $G$ be a graph that admits a tree-cut decomposition $(T,\X)$ of adhesion at most $\xi$ such that every bag has size at most $1$.
We denote $\X = (X_t: t \in V(T))$, and denote the torso at $t$ by $G_t$ for each $t \in V(T)$.
For every node $t$ of $T$, let $k_t$ be a positive integer such that $k_t+\lvert X_t \rvert \geq 2$ and $G_t-X_t$ is $k_t$-colorable with clustering $N$.

We shall prove that $G$ is $\max_{t \in V(T)}\{k_t + \lvert X_t \rvert\}$-colorable with clustering at most $N^*$. 
Let $t_1$ be a node of $T$.
We treat $T$ as a rooted tree rooted at $t_1$.
We order the nodes of $T$ as $t_1,t_2,...,t_{\lvert V(T) \rvert}$ by a depth-first-search order starting at $t_1$.

For every node $t$ of $T$, let $\downarrow t$ be the maximal subtree of $T$ rooted at $t$.
For each $i \in [\lvert V(T) \rvert]-\{1\}$, let $e_i$ be the edge of $T$ incident with $t_i$ and its parent.
For every $i \in [\lvert V(T) \rvert]-\{1\}$, we say that a subgraph $C$ of $G$ {\it crosses} $e_i$ if $V(C) \cap A_{e_i,t_i} \neq \emptyset \neq V(C) \cap B_{e_i,t_i}$.
Note that $e_1$ is undefined, but for convenience, we assume that no subgraph of $G$ crosses $e_1$.

For every induced connected subgraph $C$ of $G$, we define 
	\begin{itemize}
		\item $\gamma(C)$ to be the minimum index $j$ such that either $V(C) \cap X_{t_j} \neq \emptyset$, or $C$ contains an edge of $G$ corresponding to an edge of $G_{t_j}$, 
		\item $\tau(C) = \min\{\ell \in [\lvert V(T) \rvert]: V(C) \cap X_{t_\ell} \neq \emptyset\}$, and
		\item $\sigma(C)=(\gamma(C),\tau(C))$.
	\end{itemize}
Intuitively, $\gamma(C)$ can be considered as the ``root of $C$'' in the sense that it indicates when $C$ ``first appears''.
We compare $\sigma(C)$ by the lexicographic order.
That is, we say that $\sigma(C_1)$ is smaller than $\sigma(C_2)$, denoted by $\sigma(C_1) \preceq \sigma(C_2)$, if either $\gamma(C_1)<\gamma(C_2)$, or $\gamma(C_1)=\gamma(C_2)$ and $\tau(C_1)<\tau(C_2)$.
Since $\lvert X_t \rvert \leq 1$ for every $t \in V(T)$, $\sigma$ is a total order for any set of pairwise disjoint induced subgraphs of $G$.

A vertex $v$ of $G$ is {\it special} if $v$ is not adjacent to any vertex in $V(G)-X_{\downarrow t}$, where $t$ is the node of $T$ with $v \in X_t$.
For a node $t$ of $T$ and a vertex or subgraph $A$ of $G$, we say that $A$ is {\it $G_t$-relevant} if $A$ is incident with an edge of $G$ corresponding to an edge of $G_t$.

We shall define $\max_{t \in V(T)}\{k_t + \lvert X_t \rvert\}$-colorings $c_0,c_1^1,c_1,c_2^1,c_2,...,c_{\lvert V(T) \rvert}^1,c_{\lvert V(T) \rvert}$ of subgraphs of $G$ by a greedy algorithm such that for each $i \in [\lvert V(T) \rvert]$, 
	\begin{itemize}
		\item $c^1_{i}$ is obtained from $c_{i-1}$ by further coloring the vertex in $X_{t_i}$ (if $X_{t_i} \neq \emptyset$ and this vertex was uncolored) and all uncolored $G_{t_i}$-relevant vertices, and 
		\item $c_{i}$ is obtained from $c^1_{i}$ by further coloring all uncolored special vertices contained in $X_{\downarrow t_{i}}-X_{t_{i}}$ adjacent to some $G_{t_i}$-relevant $c_i^1$-monochromatic component. 
	\end{itemize}
Note that for each $t \in V(T)$, no special vertex in $X_{\downarrow t}-X_t$ is incident with an edge of $G_t$.

Now we formally define those colorings.
Define $c_0$ to be the coloring with empty domain.
For each $i \geq 1$, define the following.
\begin{itemize}
	\item Define $\ell$ as follows.
		\begin{itemize}
			\item If $X_{t_i} = \emptyset$, then let $\ell=0$.
			\item If $X_{t_i} \neq \emptyset$ and the vertex in $X_{t_i}$ is $c_{i-1}$-colored, then let $\ell$ be the color of this vertex.
			\item If $X_{t_i} \neq \emptyset$, the vertex in $X_{t_i}$ is $c_{i-1}$-uncolored, and no edge of $G_{t_i}-X_{t_i}$ corresponds to an edge of $G$ between a $c_{i-1}$-uncolored vertex and a $c_{i-1}$-monochromatic component crossing $e_i$, then let $\ell=1$.
			\item Otherwise, 
				\begin{itemize}
					\item let $C$ be the $c_{i-1}$-monochromatic component with minimum $\sigma(C)$ among all $c_{i-1}$-monochromatic components with the following property: $C$ crosses $e_i$ and is incident with an edge $e$ of $G$ corresponding to an edge of $G_{t_i}-X_{t_i}$ such that the other end of $e$ is $c_{i-1}$-uncolored;
					\item let $\ell$ be the color of $C$. 
				\end{itemize}
		\end{itemize}
	\item If $X_{t_i} \neq \emptyset$, then define $c^1_i(v)=\ell$, where $v$ is the vertex in $X_{t_i}$.
	\item Let $c_i'$ be a coloring of $G_{t_i}-X_{t_i}$ by using colors in $[k_{t_i}+\lvert X_{t_i} \rvert]-\{\ell\}$ with clustering $N$.
		(Note that such a coloring $c_i'$ exists since $\lvert [k_{t_i}+\lvert X_{t_i} \rvert]-\{\ell\} \rvert \geq k_{t_i}$ by the definition of $\ell$.)
	\item For each $G_{t_i}$-relevant $c_{i-1}$-uncolored vertex $v$ of $G-X_{t_i}$, 
		\begin{itemize}
			\item let $c'_i(v)$ be the color in $c_i'$ of the peripheral vertex of $G_{t_i}$ corresponding to the component of $T-t_i$ containing $v$, 
			\item if $v$ is not adjacent to any $c_{i-1}$-monochromatic component, then define $c^1_i(v)=c'_i(v)$,  
			\item otherwise, 
				\begin{itemize}
					\item let $C_v$ be the $c_{i-1}$-monochromatic component adjacent to $v$ such that $\sigma(C_v)$ is as small as possible, and
					\item define $c^1_i(v)$ to be a color in $[k_{t_i} + \lvert X_{t_i} \rvert]$ such that 
						\begin{itemize}
							\item if the color of $C_v$ is not equal to $c_i'(v)$, then $c^1_i(v)=c_i'(v)$,
							\item otherwise, $c^1_i(v)$ is an arbitrary color in $[k_{t_i} + \lvert X_{t_i} \rvert]$ different from the color of $C_v$.
						\end{itemize}
						(Note that $c^1_i(v)$ is different from the color of $C_v$ in either case.)
				\end{itemize}
		\end{itemize}
	\item For each $c^1_{i}$-uncolored special vertex $v \in X_{\downarrow t_i}-X_{t_i}$ adjacent to some $G_{t_i}$-relevant $c_i^1$-monochromatic component, 
		\begin{itemize}
			\item let $C_v$ be the $G_{t_i}$-relevant $c_i^1$-monochromatic component adjacent to $v$ such that $\sigma(C_v)$ is as small as possible, and
			\item define $c_i(v)$ to be a color in $[k_{t_i} + \lvert X_{t_i} \rvert]$ different from the color of $C_v$.
		\end{itemize}
\end{itemize}
Define $c$ to be $c_{\lvert V(T) \rvert}$.
Note that $c$ is a $\max_{t \in V(T)}\{k_t + \lvert X_t \rvert\}$-coloring of $G$.

For every $i,j \in [\lvert V(T) \rvert]$, we define the {\it $(e_i,j)$-rank} of a $c_j$-monochromatic component $C$ crossing $e_i$ to be $\alpha$ if $\sigma(C)$ is the $\alpha$-th smallest among all $c_j$-monochromatic components crossing $e_i$.

For every $i \in [\lvert V(T) \rvert]$, let $\overline{i}$ be the largest integer in $[\lvert V(T) \rvert]$ such that $t_{\overline{i}}$ is a descendant of $t_i$.

\medskip

\noindent{\bf Claim 1:} Let $i,j \in [\lvert V(T) \rvert]$ with $1 \leq i-1 \leq j$.
Let $C$ be a $c_j$-monochromatic component crossing $e_i$.
Let $j' \in [j,\overline{i}]$.
Let $C'$ be the $c_{j'}$-monochromatic component containing $C$.
Let $M'$ be a $c_{j'}$-monochromatic component crossing $e_i$.
If the $(e_i,j')$-rank of $M'$ is smaller than the $(e_i,j')$-rank of $C'$, then $M'$ contains a $c_j$-monochromatic component $M$ crossing $e_i$ such that the $(e_i,j)$-rank of $M$ is smaller than the $(e_i,j)$-rank of $C$. 

\medskip

\noindent{\bf Proof of Claim 1:}
When $j'=j$, the claim obviously holds.
So we may assume that $j'>j$.

Since $j'>j \geq i-1$, $j' \geq i$.
Since $i \leq j' \leq \overline{i}$, by the depth-search-ordering, $t_{j'}$ is a descendant of $t_i$.
Note that every vertex that is $c_{j'}$-colored but not $c_{i-1}$-colored is either $G_{t_{j''}}$-related for some $i \leq j'' \leq j'$ but not $G_{t_{j'''}}$-related for any $j''' \in [i-1]$, or a special vertex contained in $X_{\downarrow t_{j''}}$ for some $i \leq j'' \leq j'$.
So every vertex that is $c_{j'}$-colored but not $c_{i-1}$-colored is contained in $X_{\downarrow t_i}$ and is not incident with any edge in $\adh_{(T,\X)}(e_i)$.

Since $M'$ and $C$ cross $e_i$, and the $(e_i,j')$-rank of $M'$ smaller than the $(e_i,j')$-rank of $C'$, we know that $\gamma(M') \leq \gamma(C') \leq \gamma(C) \leq i-1$, and there exists a $c_{j'}$-monochromatic path $P$ contained in $M'$ containing an edge $e$ in $\adh_{(T,\X)}(e_i)$ such that either one end of $P$ is in $X_{t_{\gamma(M')}}$, or $P$ contains an edge of $G$ corresponding to an edge of $G_{t_{\gamma(M')}}$ and contains $X_{t_{\tau(M')}}$. 
We choose $e$ such that $P$ can be chosen to be as short as possible. 
So $P$ is contained in $G[e \cup (V(G)-X_{\downarrow t_i})]$.
Since $e \in \adh_{(T,\X)}(e_i)$, $e$ is an edge of the torso at the parent of $t_i$, so both ends of $e$ are $c_{i-1}$-colored and hence are $c_j$-colored.
Since every vertex that is $c_{j'}$-colored but $c_j$-uncolored is contained in $X_{\downarrow t_i}$, $P$ is a $c_j$-monochromatic path.
Then $e$ is contained in some $c_j$-monochromatic component $M$ crossing $e_i$.
Hence $M$ contains $P$. 
So $\sigma(M) \preceq \sigma(M')$. 
Since the $(e_i,j')$-rank of $M'$ is smaller than the $(e_i,j')$-rank of $C'$, we know $\sigma(M) \preceq \sigma(M') \prec \sigma(C') \preceq \sigma(C)$.
So the $(e_i,j)$-rank of $M$ is smaller than the $(e_i,j)$-rank of $C$.
Since both $M$ and $M'$ contain $P$, $M'$ contains $M$.
This proves the claim.
$\Box$

\medskip

\noindent{\bf Claim 2:} For every $i,j \in [\lvert V(T) \rvert]$ with $i \neq 1$, if $C_1$ and $C_2$ are $c_j$-monochromatic components crossing $e_i$ such that the $(e_i,j)$-rank of $C_1$ is smaller than the $(e_i,j)$-rank of $C_2$, then for every $i' \geq i$ for which $C_1$ and $C_2$ cross $e_{i'}$, the $(e_{i'},j)$-rank of $C_1$ is smaller than the $(e_{i'},j)$-rank of $C_2$.

\medskip

\noindent{\bf Proof of Claim 2:}
Let $i,j,C_1,C_2,i'$ be the ones stated in the statement of this claim.
Since the $(e_i,j)$-rank of $C_1$ is smaller than the $(e_i,j)$-rank of $C_2$, $\sigma(C_1) \prec \sigma(C_2)$. 
Since $C_1$ and $C_2$ cross $e_{i'}$, the $(e_{i'},j)$-rank of $C_1$ and $C_2$ are well-defined.
Since $\sigma(C_1) \prec \sigma(C_2)$, the $(e_{i'},j)$-rank of $C_1$ is smaller than the $(e_{i'},j)$-rank of $C_2$.
$\Box$

\medskip

\noindent{\bf Claim 3:} For every $i \in [\lvert V(T) \rvert]-\{1\}$, if $C$ is a $c_{i-1}$-monochromatic component crossing $e_i$ with $(e_i,i-1)$-rank $1$, then the $c_{\overline{i}}$-monochromatic component containing $C$ equals $C$.

\medskip

\noindent{\bf Proof of Claim 3:}
For every $\alpha \in [i-1,\overline{i}]$, let $C_\alpha$ be the $c_\alpha$-monochromatic component containing $C$.
Suppose to the contrary that $C_{\overline{i}} \neq C$.
So there exists $i^* \in [i,\overline{i}]$ such that $C_{i^*} \neq C_{i^*-1}=C$.
Hence there exists a vertex $v \in V(C_{i^*})-V(C_{i^*-1})$ adjacent to a vertex $u \in V(C_{i^*-1})$.
So $v$ is $c_{i^*-1}$-uncolored but $c_{i^*}$-colored.
Hence $v \in X_{\downarrow t_{i^*}}$.
We choose such $v$ such that $v$ is $c^1_{i^*}$-colored if possible.
Since $v$ is $c_{i^*-1}$-uncolored, $v$ is not incident with an edge in $\adh_{(T,\X)}(e_{i^*})$, so $u \in X_{\downarrow t_{i^*}}$.
Hence if $C_{i^*-1}$ does not cross $e_{i^*}$, then $V(C_{i^*-1}) \subseteq X_{\downarrow t_{i^*}}$, so $C_{i^*-1}=C$ does not cross $e_i$, a contradiction.
So $C_{i^*-1}$ crosses $e_{i^*}$.

By Claim 1, since the $(e_i,i-1)$-rank of $C$ is 1, the $(e_i,i^*-1)$-rank of $C_{i^*-1}$ is 1.
By Claim 2, the $(e_{i^*},i^*-1)$-rank of $C_{i^*-1}$ is 1.

Suppose $v \in X_{t_{i^*}}$.
Since $v$ is $c_{i^*-1}$-uncolored, $v$ is not incident with any edge in $\adh_{(T,\X)}(e_{i^*})$, so $v$ is special.
Since $i^* \geq i \geq 2$, there exists the parent $p$ of $t_{i^*}$.
Let $i_p$ be the integer such that $p=t_{i_p}$.
So $v \in X_{\downarrow p}-X_p$ is a $c_{i_p-1}$-uncolored special vertex.
Since $C_{i^*-1}$ crosses $e_i^*$, $C_{i^*-1}$ contains an edge of $G$ corresponding to an edge of $G_p$.
Since $u \in X_{\downarrow t_{i^*}}$ and $uv \in E(G)$, $u$ is not special.
Since $u$ is $c_{i^*-1}$-colored, $u$ is incident with some edge of $G$ corresponding to an edge of $G_z$ for some $z=t_{i_z}$ with $i_z \leq i^*-1$.
Since $u \in X_{\downarrow p}-X_p$, $u$ is incident with some edge of $G$ corresponding to an edge of $G_p$.
So $u$ is $c^1_{i_p}$-colored.
Since $v \in X_{\downarrow p}-X_p$ is a $c^1_{i_p}$-uncolored special vertex adjacent to $u$, and $u$ is $c^1_{i_p}$-colored and $G_p$-relevant, $v$ is $c_{i_p}$-colored so is $c_{i^*-1}$-colored, a contradiction.

Hence $v \not \in X_{t_{i^*}}$.
By the definition of $c_{i^*}$, since $C_{i^*-1}$ has $(e_{i^*},i^*-1)$-rank 1, if $v$ is $G_{t_{i^*}}$-relevant, then $c_{i^*}(v) \neq c_{i^*}(u)$, a contradiction.
So $v$ is $c_{i^*}$-colored but not $G_{t_{i^*}}$-relevant.
Hence $v$ is a special vertex in $X_{\downarrow t_{i^*}}-X_{t_{i^*}}$ and is $c^1_{i^*}$-uncolored.
Since we choose $v$ such that $v$ is $c^1_{i^*}$-colored if possible, $C_{i^*-1}$ is a $c^1_{i^*}$-monochromatic component.

Note that $C_{i^*-1}$ is the $c^1_{i^*}$-monochromatic component incident with some edge of $G_{t_{i^*}}$ crossing $e_{i^*}$ such that $\sigma(C_{i^*-1})$ is minimum. 
Note that for every $c^1_{i^*}$-monochromatic component $C$ adjacent to $v$ but not crossing $e_{i^*}$, $\gamma(C) \geq i^*> \gamma(C_{i^*-1})$, so $\sigma(C_{i^*-1}) \prec \gamma(C)$.
Hence $c_{i^*}(v) \neq c_{i^*}(u)$ by the definition of $c_{i^*}$, a contradiction.
This proves the claim.
$\Box$

\medskip

For every $i \in [\lvert V(T) \rvert]$, let $W_i$ be the set of $c_{i-1}$-colored vertices in $X_{\downarrow t_i}$.

\medskip

\noindent{\bf Claim 4:} For every $i \in [\lvert V(T) \rvert]$, $\lvert W_i \rvert \leq N_0$.

\medskip

\noindent{\bf Proof of Claim 4:}
When $i=1$, $W_i = \emptyset$.
So we may assume $i \geq 2$.
By the definition of $c_{i-1}$, every vertex in $W_i$ is either a special vertex or incident with an edge in $\adh_{(T,\X)}(e_i)$.

Let $Z$ be the set of vertices in $W_i$ incident with some edge in $\adh_{(T,\X)}(e_i)$.
So $\lvert Z \rvert \leq \lvert \adh_{(T,\X)}(e_i) \rvert \leq \xi$.

Note that every vertex in $W_i-Z$ is special.
Let $v \in W_i-Z$.
Since $v \in W_i$, there exists $i_v \in [i-1]$ such that $v$ is $c_{i_v}$-colored but not $c_{i_v-1}$-colored.
Since $v$ is not incident with any edge in $\adh_{(T,\X)}(e_i)$, $v$ is adjacent to a $c^1_{i_v}$-monochromatic component $C$.
Let $u$ be a vertex of $C$ adjacent to $v$.
Since $v$ is special, every neighbor of $v$ is not special and is contained in $X_{\downarrow t_i}$.
So $u$ is not special and is contained in $X_{\downarrow t_i}$.
Since $i_v \leq i-1$, $u \in W_i$.
Hence $u \in Z$.

Therefore, every vertex in $W_i-Z$ is adjacent to a vertex in $Z$.
For every $z \in Z$, if $z'$ is a special vertex adjacent to $z$, then the edge $zz'$ belongs to $\adh_{(T,\X)}(e_{i_z})$, where $i_z$ is the integer such that $z \in X_{t_{i_z}}$.
So there are at most $\xi \cdot \lvert Z \rvert$ special vertices adjacent to $Z$.
Hence $\lvert W_i \rvert \leq \lvert Z \rvert + \xi \lvert Z \rvert \leq \xi(\xi+1)=N_0$.
$\Box$

\medskip

\noindent{\bf Claim 5:} For every $i \in [\lvert V(T) \rvert]$, there are at most $N_1$ vertices $v$ in $X_{\downarrow t_i}$ such that $v$ is $c_i$-colored, $c_{i-1}$-uncolored, and adjacent to $W_i-X_{t_i}$. 

\medskip

\noindent{\bf Proof of Claim 5:}
Let $\Q=\{Q: Q$ is a component of $T-t_i$ with $X_Q \cap W_i \neq \emptyset\}$.
By Claim 4, $\lvert \Q \rvert \leq \lvert W_i \rvert \leq N_0$.
For each $Q \in \Q$, let $i_Q$ be the index such that $t_{i_Q}$ is the root of $Q$.

Let $Z=\{v \in X_{\downarrow t_i}: v$ is $c_i$-colored, $c_{i-1}$-uncolored, and adjacent to $W_i-X_{t_i}\}$.
For any $v \in Z-\bigcup_{Q \in \Q}X_Q$, since $v$ is adjacent to $W_i-X_{t_i}$, $v$ is incident with an edge in $\bigcup_{Q \in \Q}\adh_{(T,\X)}(e_{i_Q})$.

Let $Z' = \{v \in Z \cap \bigcup_{Q \in \Q}X_Q: v$ is $G_{t_i}$-relevant$\}$.
Let $Z'' = \{v \in Z \cap \bigcup_{Q \in \Q}X_Q: v$ is not $G_{t_i}$-relevant$\}$.
Then for every $v \in Z'$, $v$ is incident with an edge in $\bigcup_{Q \in \Q}\adh_{(T,\X)}(e_{i_Q})$.
So $\lvert Z-Z'' \rvert \leq \lvert \bigcup_{Q \in \Q}\adh_{(T,\X)}(e_{i_Q}) \rvert \leq \lvert \Q \rvert \cdot \xi \leq N_0\xi$.

Note that for every $v \in Z''$, $v$ is $c_i$-colored, $c_{i-1}$-uncolored and is not $G_{t_i}$-relevant, so $v$ is special.
Since $v \in Z$, $v$ is adjacent to $X_{\downarrow t_{i_v}} \cap W_i$, where $i_v$ is the index such that $v \in X_{t_{i_v}}$.
For every $u \in W_i$, $u$ is adjacent to at most $\lvert \adh_{(T,\X)}(e_{i_u}) \rvert \leq \xi$ vertices $u'$ with $t_{i_u} \in V(\downarrow t_{i_{u'}})$, where $i_u$ and $i_{u'}$ are the indices such that $u \in X_{t_{i_u}}$ and $u' \in X_{t_{i_{u'}}}$, respectively.
Hence $\lvert Z'' \rvert \leq \xi \cdot \lvert W_i\rvert \leq \xi N_0$ by Claim 4.
Therefore, $\lvert Z \rvert \leq \lvert Z-Z'' \rvert + \lvert Z'' \rvert \leq 2N_0\xi=N_1$.
$\Box$

\medskip

\noindent{\bf Claim 6:} For every $i \in [\lvert V(T) \rvert]$, if $C_{i-1}$ is a $c_{i-1}$-monochromatic component crossing $e_i$, and $C_i$ is the $c_i$-monochromatic component containing $C_{i-1}$, then $\lvert V(C_i)-V(C_{i-1}) \rvert \leq N_2$.

\medskip

\noindent{\bf Proof of Claim 6:}
Let $Z=V(C_i)-(V(C_{i-1}) \cup W_i)$.
Recall that $V(C_i)-V(C_{i-1}) \subseteq X_{\downarrow t_i}$.
So $\lvert V(C_i)-V(C_{i-1}) \rvert \leq \lvert Z \rvert + \lvert W_i \rvert \leq \lvert Z \rvert + N_0$ by Claim 4.
To prove this claim, it suffices to prove that $\lvert Z \rvert \leq N_2-N_0$. 

Since $V(C_i)-V(C_{i-1}) \subseteq X_{\downarrow t_i}$ and $Z \cap W_i=\emptyset$, if $v$ is in $Z$, then $v$ is $c_{i-1}$-uncolored, so all neighbors of $v$ are contained in $X_{\downarrow t_i}$.
So for every vertex $v \in Z$, one of the following holds.
	\begin{itemize}
		\item[(i)] $v \in X_{t_i}$.
		\item[(ii)] $v \not \in X_{t_i}$ is adjacent to some vertex in $W_i-X_{t_i}$.
		\item[(iii)] $v \not \in X_{t_i}$ is not adjacent to any vertex in $W_i-X_{t_i}$, and $v$ is $G_{t_i}$-relevant. 
		\item[(iv)] $v \not \in X_{t_i}$ is not adjacent to any vertex in $W_i$ and $v$ is a special vertex.
	\end{itemize}

Let $S=\{v \in X_{\downarrow t_i}: v$ is $c_i$-colored, $c_{i-1}$-uncolored, and adjacent to $W_i-X_{t_i}\}$.
Hence every vertex in $Z$ satisfying (ii) is contained in $S$.
So the number of vertices in $Z$ satisfying (i) or (ii) is at most $1+\lvert S \rvert \leq 1+N_1$ by Claim 5.

For simplicity of notations, for every $v \in X_{\downarrow t_i}-X_{t_i}$, we define $c_i'(v)$ to be $c_i'(x)$, where $x$ is the vertex in $G_{t_i}-X_{t_i}$ corresponding to the component $Q$ of $T-t_i$ with $v \in X_{Q}$.

Suppose that there exists $v \in X_{\downarrow t_i} \cap Z -W_i$ satisfying (iii) such that $c_i(v) \neq c_i'(v)$. 
Then $v$ is adjacent to some $c_{i-1}$-monochromatic component by the definition of $c_i^1$ and $c_i$.
Since all neighbors of $v$ are contained in $X_{\downarrow t_i}$, $v$ is adjacent to some vertex in $W_i$.
Since $v$ satisfies (iii), $\emptyset \neq X_{t_i} \subseteq W_i$, and the $c_{i-1}$-monochromatic component containing $X_{t_i}$ is the only $c_{i-1}$-monochromatic component adjacent to $v$.
Since the color of $X_{t_i}$ is not in the image of $c_i'$, $c_i(v)=c_i^1(v)=c_i'(v)$, a contradiction. 

This shows that 
	\begin{itemize}
		\item[(a)] for every $v \in X_{\downarrow t_i} \cap Z-W_i$ satisfying (iii), $c_i(v)=c_i'(v)$, so either $X_{t_i}=\emptyset$, or $c_i(v)=c_i'(v)$ is different from the color of the vertex in $X_{t_i}$.
	\end{itemize}

We first assume that $X_{t_i} \neq \emptyset$ and the color of the vertex in $X_{t_i}$ equals the color of $C_i$.
So no vertex in $Z$ satisfies (iii) by (a).
Let $z$ be a vertex in $Z$ satisfying (iv).
Since $z \in V(C_i)-V(C_{i-1})$ and $z$ is $c_{i-1}$-uncolored, $z$ is adjacent to a vertex $u_z$ in $V(C_i)$ with $c_i(z)=c_i(u_z)$.
Since $z$ satisfies (iv), $u_z \not \in W_i$, so $u_z \in Z$.
Since $z$ is special, $u_z$ satisfies (ii) or (iii).
Since $u_z \in Z$ and no vertex in $Z$ satisfies (iii), $u_z$ satisfies (ii), so $u_z \in S$.
Note that $u_z$ is adjacent to at most $\xi$ special vertices.
Hence the number of vertices in $Z$ satisfying (iv) is at most $\lvert S \rvert \cdot \xi \leq N_1 \xi$ by Claim 5.
Therefore $\lvert Z \rvert \leq 1+N_1+0+N_1 \xi \leq N_2-N_0$ and we are done.

Hence we may assume that either $X_{t_i}=\emptyset$, or $X_{t_i} \neq \emptyset$ and the color of the vertex in $X_{t_i}$ is different from the color of $C_i$.
So $V(C_{i-1}) \cap X_{\downarrow t_i} \subseteq W_i-X_{t_i}$ and $V(C_i) \cap X_{t_i} = \emptyset$.

We define a {\it $Z$-component} to be a component of $C_i-(V(C_{i-1}) \cup W_i)$.
That is, every $Z$-component is a component of $G[Z]$.
Let $Z_2 = \{v \in Z: v$ satisfies (ii)$\}$.

Note that every vertex in $X_{\downarrow t_i}$ incident with some edge in $\adh_{(T,\X)}(e_i)$ is $c_{i-1}$-colored.
Since $V(C_{i-1}) \cap X_{\downarrow t_i} \subseteq W_i-X_{t_i}$ and $C_i$ is connected, for every vertex $v$ in $V(C_i)-V(C_{i-1})$, there exists a path in $C_i$ from $v$ to $W_i-X_{t_i}$ internally disjoint from $W_i-X_{t_i}$.
So every $Z$-component contains a vertex in $Z_2$.

Hence 
	\begin{itemize}
		\item[(b)] there are at most $\lvert Z_2 \rvert \leq \lvert S \rvert \leq N_1$ $Z$-components.
	\end{itemize}

Let $Z_{3,4}=\{v \in Z: v$ satisfies (iii) or (iv)$\}$.
We define a {\it $Z_{3,4}$-component} to be a component of $G[Z_{3,4}]$.

Let $v \in Z_2$.
Let $i_v$ be the index such that $v \in X_{t_{i_v}}$.
Then $v$ is adjacent to at most $\xi$ vertices in $V(G)-X_{\downarrow t_{i_v}}$.
Let $u$ be a neighbor of $v$ contained in $Z_{3,4} \cap X_{\downarrow t_{i_v}}$. 
Then $u$ is not special.
So $u$ satisfies (iii).
Hence $u$ is $G_{t_i}$-relevant, so $u$ is incident with an edge in $\adh_{(T,\X)}(e_{i_v})$.

Therefore, every vertex in $Z_2$ is adjacent to at most $2\xi$ vertices in $Z_{3,4}$.
So 
	\begin{itemize}
		\item [(c)] every $Z$-component consists of at most $\lvert Z_2 \rvert$ vertices in $Z_2$ and at most $\lvert Z_2 \rvert \cdot 2\xi$ $Z_{3,4}$-components.
	\end{itemize}

Let $Z_3 = \{z \in Z: z$ satisfies (iii)$\}$.
We define a {\it $Z_3$-component} to be a component of $G[Z_3]$.
Let $Z_4 = \{z \in Z: z$ satisfies (iv)$\}$.
Let $Z_4' = \{z \in Z_4: z$ is adjacent to a $Z_3$-component$\}$.
Since every vertex in $Z_4$ is special, vertices in $Z_4$ are pairwise nonadjacent.
So every $Z_{3,4}$-component intersecting $Z_4-Z_4'$ consists of one vertex in $Z_4-Z_4'$.
We say that a $Z_{3,4}$-component is a {\it $Z_{3,4}'$-component} if it is disjoint from $Z_4-Z_4'$.

If $v$ is a vertex in $Z_4-Z_4'$, then since $v$ is special and not adjacent to any vertex in $W_i$, $v$ is adjacent to $Z_2$, for otherwise the $Z$-component containing $v$ consists of $v$ and hence does not contain a vertex satisfying (ii), a contradiction.
Since every vertex in $Z_4-Z_4'$ is special, $\lvert Z_4 - Z_4' \rvert \leq \lvert Z_2 \rvert \xi$.

For any $Z_3$-component $Q$, we define a {\it $Q$-branch} to be a component $Q'$ of $\downarrow t_i-t_i$ such that $X_{Q'} \cap V(Q) \neq \emptyset$.
Since every vertex in $Z_4'$ is adjacent to a $Z_3$-component and every vertex in $Z_4$ is special, every component $Q'$ of $\downarrow t_i-t_i$ with $X_{Q'} \cap Z_4' \neq \emptyset$ is a $Q$-branch for some $Z_3$-component $Q$.
Since every vertex in $Z_3$ is $G_{t_i}$-relevant, for every component $Q'$ of $\downarrow t_i-t_i$, $X_{Q'}$ contains at most $\xi$ vertices in $Z_3$.
Since every vertex in $Z_4'$ is special, every vertex is adjacent to at most $\xi$ vertices in $Z_4'$.

Similarly, for any $Z'_{3,4}$-component $Q$, we define a {\it $Q$-branch} to be a component $Q'$ of $\downarrow t_i-t_i$ such that $X_{Q'} \cap V(Q) \neq \emptyset$.
For every $Z'_{3,4}$-component $Q$ and every $Q$-branch $Q'$, $Q'$ is also a $Q_3$-branch for some $Z_3$-component $Q_3$ by the definition of $Z_4'$, since every vertex in $Z_4$ is special.
Hence, for every $Q$-branch $Q'$ for some $Z'_{3,4}$-component $Q$, as shown in the previous paragraph, $X_{Q'} \cap (Z_3 \cup Z_4')$ consists of at most $\xi$ vertices in $Z_3$ and at most $\lvert X_{Q'} \cap Z_3 \rvert \cdot \xi \leq \xi^2$ vertices in $Z_4'$.
So 
	\begin{itemize}
		\item[(d)] for every $Z'_{3,4}$-component $Q$, $\lvert V(Q) \rvert$ is at most $\xi^2+\xi$ times the number of $Q$-branches.
	\end{itemize}

Let $Q$ be a $Z'_{3,4}$-component.
Note that $V(Q) \cap X_{t_i} = \emptyset$ as they have different colors.
Let $Q'$ be the subgraph of $G_{t_i}-X_{t_i}$ induced by the vertices of $G_{t_i}-X_{t_i}$ corresponding to the $Q$-branches. 
Note that $Q'$ can be obtained from $Q$ by identifying vertices.
Since $Q$ is connected, $Q'$ is connected.
By (a), every vertex in $Z$ satisfying (iii) satisfies $c_i(v)=c_i'(v)$, so $Q'$ is contained in a $c_i'$-monochromatic component in $G_{t_i}-X_{t_i}$.
Hence $Q'$ contains at most $N$ vertices by the definition of $c_i'$.
So there are at most $N$ $Q$-branches. 
By (d), $\lvert V(Q) \rvert \leq (\xi^2+\xi) \cdot N$.

Since every $Z_{3,4}$-component either consists of one vertex in $Z_4-Z_4'$ or is a $Z_{3,4}'$-component, every $Z_{3,4}$-component contains at most $(\xi^2+\xi)N$ vertices.
So by (c), every $Z$-component contains at most $\lvert Z_2 \rvert + \lvert Z_2 \rvert \cdot 2\xi \cdot (\xi^2+\xi)N \leq \lvert S \rvert \cdot (1+2\xi^2(\xi+1)N) \leq N_1 \cdot (1+2\xi^2(\xi+1)N)$ vertices.
By (b), $\lvert Z \rvert \leq N_1 \cdot (1+2\xi^2(\xi+1)N) \cdot N_1 = (1+2\xi^2(\xi+1)N)N_1^2 \leq N_2-N_0$.
This proves the claim.
$\Box$

\medskip

\noindent{\bf Claim 7:} For every $i \in [\lvert V(T) \rvert]-\{1\}$ and $k \in {\mathbb N}$, if $C$ is a $c_{i-1}$-monochromatic component crossing $e_i$ and with $(e_i,i-1)$-rank $k$, then $\lvert V(M) \cap X_{\downarrow t_i} \rvert \leq f(k)$, where $M$ is the $c$-monochromatic component containing $C$. 

\medskip

\noindent{\bf Proof of Claim 7:}
Let $i,k,C,M$ be the ones as stated in the claim.
We shall prove this claim by induction on $k$.
For every $\alpha \in [i-1,\overline{i}]$, let $C_\alpha$ be the $c_\alpha$-monochromatic component containing $C$.
Note that $C=C_{i-1}$ and $\lvert V(M) \cap X_{\downarrow t_i} \rvert = \lvert V(C_{\overline{i}}) \cap X_{\downarrow t_i} \rvert$.

If $k=1$, then by Claims 3 and 4, $\lvert V(M) \cap X_{\downarrow t_i} \rvert = \lvert V(C_{\overline{i}}) \cap X_{\downarrow t_i} \rvert = \lvert V(C) \cap X_{\downarrow t_i} \rvert \leq \lvert W_i \rvert \leq N_0=f(1)$.

So we may assume that $k \geq 2$, and the claim holds if $k$ is smaller.
For every $\alpha \in [k]$, let $M_\alpha$ be the $c$-monochromatic component containing the $c_{i-1}$-monochromatic component crossing $e_i$ with $(e_i,i-1)$-rank $\alpha$.
Let $J = \{j \in [i,\overline{i}]: C_j \neq C_{j-1}\}$.

Let $j \in J$.
Let $S_j=\{v \in V(C_j) \cap X_{\downarrow t_i}-V(C_{j-1}): v$ is adjacent to some vertex in $C_{j-1}\}$.

Let $v \in S_j$.
So there exists $u \in V(C_{j-1})$ such that $uv \in E(G)$.
Note that $v$ is $c_j$-colored but not $c_{j-1}$-colored.
So either $v \in X_{t_j}$ and $v$ is special, or $v \in X_{\downarrow t_j}-X_{t_j}$.

Suppose that $v$ is special.
Then $u$ is not special and $u \in X_{\downarrow t_v}$, where $t_v$ is the node of $T$ with $v \in X_{t_v}$.
Let $i_u$ be the index such that $u$ is $c_{i_u}$-colored but $c_{i_{u-1}}$-uncolored.
So $i_u \leq j-1$.
Since $u$ is not special, $u$ is $G_{t_{i_u}}$-relevant and $c^1_{i_u}$-colored.
Since $u$ is $c_{i_{u}}$-colored but $c_{i_{u}-1}$-uncolored, $u \in X_{\downarrow t_{i_u}}$.
Since $i_u \leq j-1$, $v \in X_{\downarrow t_{i_u}}-X_{t_{i_u}}$.
So $v$ is $c_{i_u}$-colored by the definition of $c_{i_u}$.
Hence $v$ is $c_{j-1}$-colored, a contradiction.

So $v$ is not special and $v \in X_{\downarrow t_j}-X_{t_j}$.
Since $v$ is $c_{j-1}$-uncolored, all neighbors of $v$ are contained in $X_{\downarrow t_j}$, so $u \in X_{\downarrow t_j}$.
Since $C_{j-1}$ contains $u$ and $C$, and $C$ crosses $e_i$, we know that $C_{j-1}$ crosses $e_j$.
Since $v$ is not special and is $c_j$-colored, $v$ is $c_j^1$-colored.
Since $u$ is adjacent to $v$ and is contained in $C_{j-1}$, and $c_j(u)=c_j(v)$, we know that $v$ is adjacent to a vertex $x$ of $G$ contained in a $c_{j-1}$-monochromatic component $C_x$ such that $\sigma(C_x) \prec \sigma(C_{j-1})$, by the definition of $c^1_j(v)$. 
Since $C$ crosses $e_i$, $C_{j-1}$ crosses $e_i$.
Since $\sigma(C_x) \prec \sigma(C_{j-1})$ and $x \in X_{\downarrow t_j}$, $C_x$ crosses $e_i$ and the $(e_i,j-1)$-rank of $C_x$ is smaller than the $(e_i,j-1)$-rank of $C_{j-1}$.
By Claim 1, $C_x$ contains a $c_{i-1}$-monochromatic component $C_x'$ crossing $e_i$ such that the $(e_i,i-1)$-rank of $C_x'$ is smaller than the $(e_i,i-1)$-rank of $C$.
So there exists $\alpha_x \in [k-1]$ such that $M_{\alpha_x}$ contains $C'_x$.
Since $C_x$ contains $x$ and $C_x'$, and $M_{\alpha_x}$ is monochromatic, $M_{\alpha_x}$ contains $C_x$ and $x$.

Let $i_x$ be the index such that $x \in X_{t_{i_x}}$.
Since all neighbors of $v$ are contained in $X_{\downarrow t_j}$, $t_{i_x} \in V(\downarrow t_j)$. 

Suppose $i_x \neq j$ and $v$ is not incident with an edge in $\adh_{(T,\X)}(e_{i_x})$.
So $x \in X_{\downarrow t_j}- X_{t_j}$ and $v \in X_{\downarrow t_{i_x}}-X_{t_{i_x}}$, where $t_v$ is the node of $T$ with $v \in X_{t_v}$.
Since $v$ is $c_j$-colored but not special, and since $v \in X_{\downarrow t_j}-X_{t_j}$, $v$ is $G_{t_j}$-relevant.
Since $t_{i_x} \in V(\downarrow t_j)-\{t_j\}$ and $v \in X_{\downarrow t_{i_x}}-X_{t_{i_x}}$, $v$ is incident with an edge in $\adh_{(T,\X)}(e_{i_x})$, a contradiction.

Hence either $i_x=j$ or $v$ is incident with an edge in $\adh_{(T,\X)}(e_{i_x})$.
Note that $x \in X_{\downarrow t_j} \subseteq X_{\downarrow t_i}$.

This shows that for every $\alpha \in J$, there exists $x_\alpha \in \bigcup_{\beta=1}^{k-1}(V(M_\beta) \cap X_{\downarrow t_i})$ such that either $\alpha = i_{x_\alpha}$ or there exists an edge of $G$ in $\adh_{(T,\X)}(e_{i_{x_\alpha}})$ incident with a vertex in $V(C_{\alpha})-V(C_{\alpha-1})$, where $i_{x_\alpha}$ is the index such that $x_\alpha \in X_{t_{i_{x_\alpha}}}$.
Therefore, $\lvert J \rvert \leq \lvert \bigcup_{\beta=1}^{k-1}(V(M_\beta) \cap X_{\downarrow t_i}) \rvert \cdot (\xi+1) \leq (\xi+1) \cdot \sum_{\beta=1}^{k-1}f(\beta)$ by the induction hypothesis.

Note that for every $\alpha \in J$, since $C$ crosses $e_i$ and $C_\alpha \neq C_{\alpha-1}$, we know that $C_{\alpha-1}$ crosses $e_\alpha$.
By Claim 6, for every $\alpha \in J$, $\lvert V(C_\alpha)-V(C_{\alpha-1}) \rvert \leq N_2$.
Therefore, $\lvert V(C_{\overline{i}})-V(C) \rvert = \lvert V(C_{\overline{i}})-V(C_{i-1}) \rvert = \sum_{\alpha \in J} \lvert V(C_\alpha)-V(C_{\alpha-1}) \rvert \leq \lvert J \rvert \cdot N_2 \leq (\xi+1) N_2 \cdot \sum_{\beta=1}^{k-1}f(\beta) = f(k)$.
$\Box$

\medskip

Let $M$ be a $c$-monochromatic component.
To prove this lemma, it suffices to show that $\lvert V(M) \rvert \leq N^*$.

Let $I$ be the minimal set of nodes of $T$ such that for every $t \in I$, $X_t \cap V(M) \neq \emptyset$, and for every $t' \in V(T)$ with $V(M) \cap X_{t'} \neq \emptyset$, $t'$ is a descendant of some node in $I$.
Let $c(M)$ be the color of $M$ in $c$.

\medskip

\noindent{\bf Claim 8:} If $\lvert I \rvert =1$, then $\lvert V(M) \rvert \leq N^*$.

\medskip

\noindent{\bf Proof of Claim 8:}
Let $t$ be the node in $I$, and let $i$ be the index such that $t=t_i$.
Let $\Q = \{Q: Q$ is a component of $\downarrow t_i-t_i$ with $X_{Q} \cap V(M) \neq \emptyset\}$.
We may assume $V(M) \neq X_{t_i}$, for otherwise $\lvert V(M) \rvert = \lvert X_{t_i} \rvert=1 \leq N^*$ and we are done.
So $\Q \neq \emptyset$.
For each $Q \in \Q$, let $i_Q$ be the index such that $Q = \downarrow t_{i_Q}$.
Note that $\lvert V(M) \rvert = \lvert X_{t_i} \rvert + \sum_{Q \in \Q} \lvert V(M) \cap X_{\downarrow t_{i_Q}} \rvert$.

For every $Q \in \Q$, there exists a path in $M$ from $X_{t_i}$ to $X_{Q} \cap V(M)$, so some vertex $v_Q$ in $V(M) \cap X_{Q}$ is $G_{t_i}$-relevant.
So for each $Q \in \Q$, there exists a $c_{i_Q-1}$-monochromatic component $C_Q$ crossing $e_{i_Q}$ contained in $M$.
Note that for each $Q \in \Q$, $M$ is the $c$-monochromatic component containing $C_Q$, and the $(e_{i_Q},i_Q-1)$-rank of $c_Q$ is at most $\xi$, so $\lvert V(M) \cap X_{\downarrow t_{i_Q}} \rvert \leq f(\xi)$ by Claim 7.
Hence $\lvert V(M) \rvert = \lvert X_{t_i} \rvert + \sum_{Q \in \Q} \lvert V(M) \cap X_{\downarrow t_{i_Q}} \rvert \leq 1+\lvert \Q \rvert \cdot f(\xi)$.

Let $\W = \{W: W$ is a component of $\downarrow t_i-t_i$ with $X_{W} \cap W_i \neq \emptyset\}$.
For each $W \in \W$, let $e_W$ be the edge of $T$ between $W$ and $t_i$.

Let $Q \in \Q-\W$.
Since $v_Q$ is $G_{t_i}$-relevant, $v_Q$ is $c^1_i$-colored.
Since $Q \not \in \W$, $v_Q$ is $c_{i-1}$-uncolored.
Since $V(M) \cap X_{t_i} \neq \emptyset$ and $c(M)=c(v_Q)$, $c_i(v_Q) \neq c_i'(v_Q)$.
So $v_Q$ is adjacent to a vertex $u_Q$ in a $c_{i-1}$-monochromatic component whose color is different from $c(M)$.
Since $v_Q$ is $c_{i-1}$-uncolored, $u_Q \in X_{\downarrow t_i}$.
Since $c(u_Q) \neq c(M)$, $u_Q \in X_{\downarrow t_i}-X_{t_i}$.
So $u_Q \in \bigcup_{W \in \W}X_{W}$.
Since $Q \in \Q-\W$, $v_Q$ is incident with an edge in $\bigcup_{W \in \W}\adh_{(T,\X)}(e_W)$.

Therefore, $\lvert \Q-\W \rvert \leq \lvert \bigcup_{W \in \W}\adh_{(T,\X)}(e_W) \rvert \leq \xi \lvert \W \rvert$.
So $\lvert \Q \rvert \leq \lvert \W \rvert + \xi\lvert \W \rvert \leq (\xi+1)\lvert W_i \rvert \leq (\xi+1)N_0$ by Claim 4.
Hence $\lvert V(M) \rvert \leq 1+\lvert \Q \rvert \cdot f(\xi) \leq 1+(\xi+1)N_0 \cdot f(\xi) \leq N^*$.
$\Box$

\medskip

\noindent{\bf Claim 9:} If $\lvert I \rvert \geq 2$, then $\lvert V(M) \rvert \leq N^*$.

\medskip

\noindent{\bf Proof of Claim 9:}
Let $i$ be the largest integer such that $t_i$ is an ancestor of all nodes in $I$.
Note that such $i$ exists since $1$ is a candidate.
Let $\Q = \{Q: Q$ is a component of $\downarrow t_i-t_i$ with $X_{Q} \cap V(M) \neq \emptyset\}$.
Since $\lvert I \rvert \geq 2$, $V(M) \cap X_{t_i} = \emptyset$ and $\Q \neq \emptyset$.
By the maximality of $i$, $\lvert \Q \rvert \geq 2$.
For each $Q \in \Q$, let $i_Q$ be the index such that $Q = \downarrow t_{i_Q}$.
Note that $\lvert V(M) \rvert = \sum_{Q \in \Q} \lvert V(M) \cap X_{\downarrow t_{i_Q}} \rvert$. 

Let $Q^* \in \Q$.
For every $Q \in \Q-\{Q^*\}$, there exists a path in $M$ from $X_{Q^*} \cap V(M)$ to $X_{Q} \cap V(M)$, so some vertex $v_Q$ in $V(M) \cap X_{Q}$ is incident with an edge of $G$ corresponding to an edge of $G_{t_i}$.
Similarly, some vertex $v_{Q^*}$ in $X_{Q^*} \cap V(M)$ is incident with an edge of $G$ corresponding to an edge of $G_{t_i}$.
So for each $Q \in \Q$, there exists a $c_{i_Q-1}$-monochromatic component $C_Q$ crossing $e_{i_Q}$ contained in $M$.
Note that for each $Q \in \Q$, $M$ is the $c$-monochromatic component containing $C_Q$, so $\lvert V(M) \cap X_{\downarrow t_{i_Q}} \rvert \leq f(\xi)$ by Claim 7.
Hence $\lvert V(M) \rvert = \sum_{Q \in \Q} \lvert V(M) \cap X_{\downarrow t_{i_Q}} \rvert \leq \lvert \Q \rvert \cdot f(\xi)$. 

For each $Q \in \Q$, let $c_i'(Q)=c_i'(x_Q)$, where $x_Q$ is the peripheral vertex in $G_{t_i}-X_{t_i}$ corresponding to $Q$.
Let $\Q_1 = \{Q \in \Q: c_i'(Q) \neq c(M)\}$.
Let $\Q_2 = \Q-\Q_1$.
Let $\W = \{W: W$ is a component of $\downarrow t_i-t_i$ with $X_{W} \cap W_i \neq \emptyset\}$.
For each $W \in \W \cup \Q$, let $e_W$ be the edge of $T$ between $W$ and $t_i$.

Let $Q \in \Q_1-\W$.
Since $v_Q$ is incident with an edge of $G$ corresponding to an edge of $G_{t_i}$, $v_Q$ is $G_{t_i}$-relevant and hence is $c^1_i$-colored.
Since $Q \not \in \W$, $v_Q$ is $c_{i-1}$-uncolored.
Since $Q \in \Q_1$, $c_i(v_Q)=c(M) \neq c_i'(Q)$.
So $v_Q$ is adjacent to a vertex $u_Q$ in a $c_{i-1}$-monochromatic component whose color is $c_i'(Q)$.
Since $v_Q$ is $c_{i-1}$-uncolored, $u_Q \in X_{\downarrow t_i}$.
Since $c(u_Q) = c_i'(Q)$, $u_Q \in X_{\downarrow t_i}-X_{t_i}$ by the definition of $c_i'$.
So $u_Q \in \bigcup_{W \in \W}X_{W}$.
Since $v_Q$ is $G_{t_i}$-relevant, $v_Q$ is incident with an edge in $\bigcup_{W \in \W}\adh_{(T,\X)}(e_W)$.

Therefore, $\lvert \Q_1-\W \rvert \leq \lvert \bigcup_{W \in \W}\adh_{(T,\X)}(e_W) \rvert \leq \xi \lvert \W \rvert$.
So $\lvert \Q_1 \rvert \leq (\xi+1)\lvert \W \rvert \leq (\xi+1)\lvert W_i \rvert \leq (\xi+1)N_0$ by Claim 4.

For $\C \in \{\Q,\Q_1,\Q_2\}$, let $M_\C$ be the graph obtained from $M[V(M) \cap \bigcup_{Q \in \C}X_Q]$ by for each $Q \in \C$, identifying $Q$ into a single vertex.
Note that $M_\Q$, $M_{\Q_1}$ and $M_{\Q_2}$ are subgraphs of $G_{t_i}-X_{t_i}$.
Since $M$ is connected, $M_\Q$ is connected.
Note that $\lvert V(M_{\Q_2}) \rvert = \lvert \Q_2 \rvert$ and $M_{\Q_2} = M_\Q - V(M_{\Q_1})$.
In addition, $M_{\Q_2}$ is contained in a (not necessarily connected) $c_i'$-monochromatic subgraph of $G_{t_i}-X_{t_i}$.

We claim that $\lvert V(M_{\Q_2}) \rvert \leq \xi(2\xi+1)N_0N$.
If $\Q_1 = \emptyset$, then $M_{\Q_2}=M_\Q$ is a $c_i'$-monochromatic component of $G_{t_i}-X_{t_i}$, so $\lvert V(M_{\Q_2}) \rvert \leq N$ by the definition of $c_i'$.
So we may assume that $\Q_1 \neq \emptyset$.
For every component $R$ of $M_{\Q_2}$, let $\Q_R$ be the subset of $\Q_2$ consisting of the members of $\Q_2$ corresponding to vertices of $R$.
Since for every component $R$ of $M_{\Q_2}$, there exists a path in $M_\Q$ from $R$ to $M_{\Q_1}$ internally disjoint from $V(M_{\Q_1})$, so there exists an edge $e_R$ of $G$  between $\bigcup_{Q_R \in \Q_R}X_{Q_R}$ and $\bigcup_{Q_1 \in \Q_1}X_{Q_1}$. 
Note that $e_R \in \bigcup_{Q_1 \in \Q_1}\adh_{(T,\X)}(e_{Q_1})$.
So the number of components of $M_{\Q_2}$ is at most $\lvert \bigcup_{Q_1 \in \Q_1}\adh_{(T,\X)}(e_{Q_1}) \rvert \leq \lvert \Q_1 \rvert \xi \leq \xi(\xi+1)N_0$.
Since each component of $M_{\Q_2}$ is contained in a $c_i'$-monochromatic component of $G_{t_i}-X_{t_i}$, it contains at most $N$ vertices by the definition of $c_i'$.
Therefore, $\lvert V(M_{\Q_2}) \rvert \leq N \cdot \xi(\xi+1)N_0$.

Hence $\lvert \Q_2 \rvert = \lvert V(M_{\Q_2}) \rvert \leq \xi(\xi+1)N_0N$.
So $\lvert \Q \rvert = \lvert \Q_1 \rvert + \lvert \Q_2 \rvert \leq (\xi+1)N_0+\xi(\xi+1)N_0N = (1+N\xi)(\xi+1)N_0$.
Therefore, $\lvert V(M) \rvert \leq \lvert \Q \rvert \cdot f(\xi) \leq (1+N\xi)(\xi+1)N_0 \cdot f(\xi) \leq N^*$. 
$\Box$

\medskip

By Claims 8 and 9, $\lvert V(M) \rvert \leq N^*$.
This proves the lemma.
\end{pf}

\bigskip

Now we are ready to prove Theorem \ref{tree_cut_intro}.
The following is a restatement.

\begin{lemma} \label{coloring_const_bag}
For any positive integers $N,\xi$ and $\alpha$, there exists a positive integer $N^*=N^*(N,\xi,\alpha)$ such that the following holds.
Let $G$ be a graph that admits a tree-cut decomposition $(T,\X=(X_t: t \in V(T)))$ of adhesion at most $\xi$ such that every bag contains at most $\alpha$ vertices.
For every $t \in V(T)$, let $k_t$ be a positive integer with $k_t+\lvert X_t \rvert \geq 2$ such that the graph obtained from the torso at $t$ by deleting $X_t$ is $k_t$-colorable with clustering $N$.
Then $G$ is $\max_{t \in V(T)}\{k_t + \min\{\lvert X_t \rvert,1\}\}$-colorable with clustering $N^*$. 
\end{lemma}

\begin{pf}
Let $N,\xi$ and $\alpha$ be positive integers.
Let $\beta = N_{\ref{coloring_bag1}}(N,\xi)$, where $N_{\ref{coloring_bag1}}$ is the number $N^*$ mentioned in Lemma \ref{coloring_bag1}.
Define $N^* = \alpha \beta$.

Let $G$ be a graph that admits a tree-cut decomposition $(T,\X)$ of adhesion at most $\xi$ such that every bag contains at most $\alpha$ vertices.
Denote $\X$ by $(X_t: t \in V(T))$.
For every $t \in V(T)$, let $G_t$ be the torso at $t$ in $(T,\X)$.
Let $k_t$ be a positive integer with $k_t+\lvert X_t \rvert \geq 2$ such that $G_t-X_t$ is $k_t$-colorable with clustering $N$.

Let $G'$ be the graph obtained from $G$ by for each $t \in V(T)$ with $X_t \neq \emptyset$, identifying all vertices in $X_t$ into a vertex $v_t$.
For every $t \in V(T)$ with $X_t \neq \emptyset$, let $X'_t = \{v_t\}$; for every $t \in V(T)$ with $X_t=\emptyset$, let $X'_t=\emptyset$.
Let $\X' = (X_t': t\in V(T))$.
Then $(T,\X')$ is a tree-cut decomposition of $G'$ of adhesion at most $\xi$ such that every bag contains at most 1 vertex.
For every $t \in V(T)$, let $G'_t$ be the torso at $t$ in $(T,\X')$.
Note that for every $t \in V(T)$, $G_t-X_t = G'_t-X'_t$, so $G'_t-X'_t$ is $k_t$-colorable with clustering $N$.
If there exists $t \in V(T)$ with $k_t+\lvert X'_t \rvert \leq 1$, then since $k_t$ is a positive integer, $X_t' = \emptyset$, so $X_t = \emptyset$ and hence $k_t+\lvert X_t \rvert = k_t+\lvert X'_t \rvert \leq 1$, a contradiction.
Hence for every $t \in V(T)$, $k_t+\lvert X'_t \rvert \geq 2$.
Then by Lemma \ref{coloring_bag1}, there exists a $\max_{t \in V(T)}\{k_t + \lvert X'_t \rvert\}$-coloring $c'$ of $G'$ with clustering $\beta$. 

Define $c$ to be a coloring of $G$ such that for every $u \in V(G)$, define $c(u)=c'(v_{t_u})$, where $t_u$ is the node of $T$ with $u \in X_{t_u}$.
So $c$ is a $\max_{t \in V(T)}\{k_t + \lvert X'_t \rvert\}$-coloring of $G$.
Note that for each $t \in V(T)$, $\lvert X'_t \rvert \in \{0,1\}$, so $k_t +\lvert X'_t \rvert = k_t + \min\{\lvert X'_t \rvert,1\} = k_t + \min\{\lvert X_t \rvert,1\}$.
Hence $c$ is a $\max_{t \in V(T)}\{k_t + \min\{\lvert X_t \rvert,1\}\}$-coloring of $G$.

Let $C$ be a $c$-monochromatic component.
Let $C'$ be the graph obtained from $C$ by for each $t \in V(T)$ with $X_t \neq \emptyset$, identifying $X_t \cap V(C)$ into a vertex.
Then $C'$ is contained in a $c'$-monochromatic component in $G'$.
So $\lvert V(C') \rvert \leq \beta$.
Hence $\lvert V(C) \rvert \leq \beta \cdot \max_{t \in V(T)}\lvert X_t \rvert \leq \beta \alpha=N^*$.
So $c$ is a $\max_{t \in V(T)}\{k_t + \min\{\lvert X_t \rvert,1\}\}$-coloring of $G$ with clustering $N^*$.
\end{pf}

\section{Proof of Theorem \ref{main_intro}} \label{sec:proof}

In this section we will prove the main theorem for clustered coloring immersion-free graphs.
We need the following structure theorem proved in \cite{l_imm_str}.
A graph is {\it exceptional} if it contains exactly one vertex of degree at least two, and this vertex is incident with a loop.

\begin{theorem}[{\cite[Theorem 4.6]{l_imm_str}}] \label{global}
For any positive integers $d,h$, there exist integers $\eta=\eta(d,h)$ and $\xi=\xi(d,h)$ such that the following holds.
Let $H$ be a graph on $h$ vertices with maximum degree $d$. 
Let $G$ be a graph with no edge-cut of order exactly 3 such that $G$ does not contain an $H$-immersion.
Define $H'=H$ if $H$ is non-exceptional; otherwise, define $H'$ to be a graph obtained from $H$ by subdividing one edge.
Then there exists a tree-cut decomposition $(T,\X)$ of $G$ of adhesion at most $\eta$ such that for every $t \in V(T)$, there exists $Z_t \subseteq E(G)$ with $\lvert Z_t \rvert \leq \xi$ such that if $G_t$ is the torso at $t$, then there exists a nonnegative integer $d_t \leq d$ such that 
	\begin{enumerate}
		\item the number of vertices of degree at least $d_t$ in $G_t-Z_t$ is less than the number of vertices of degree at least $d_t$ in $H'$, 
		\item every vertex of $G_t$ of degree at least $d_t$ in $G_t-Z_t$ is a non-peripheral vertex of $G_t$,
		\item if $\lvert V(T) \rvert=1$ or $t$ is not a leaf, then every vertex in $X_t$ has degree at least $d_t$ in $G_t-Z_t$, and
		\item if $t$ is a leaf and $\lvert V(T) \rvert \geq 2$, then $\lvert X_t \rvert \leq 1$.
	\end{enumerate}
\end{theorem}

Before proving the main theorem, we need the following easy lemma.

\begin{lemma} \label{coloring_adding_edge}
Let $G$ be a graph.
Let $\xi,k,N$ be positive integers.
Let $G'$ be a graph that can be obtained from $G$ by deleting at most $\xi$ edges.
If $G'$ is $k$-colorable with clustering $N$, then $G$ is $k$-colorable with clustering $(\xi+1) N$.
\end{lemma}

\begin{pf}
Let $c$ be a $k$-coloring of $G'$ of clustering $N$.
Let $Z \subseteq E(G)$ with $\lvert Z \rvert \leq \xi$ such that $G'=G-Z$.
Since $V(G)=V(G')$, $c$ is a $k$-coloring of $G$.
If $M$ is a $c$-monochromatic component in $G$, then $M-Z$ intersects at most $\xi+1$ $c$-monochromatic components of $G'$, so $\lvert V(M) \rvert \leq (\xi+1) N$.
\end{pf}

\bigskip

Recall that $\chi_*: {\mathbb N} \cup \{0\} \rightarrow {\mathbb N}$ is the function such that for every $x \in {\mathbb N} \cup \{0\}$, $\chi_*(x)$ is the minimum $k$ such that there exists $N_x \in {\mathbb N}$ such that every graph of maximum degree at most $x$ is $k$-colorable with clustering $N_x$.

\begin{lemma} \label{coloring_upper_bdd}
For any positive integers $d \geq 3,h$, there exists a positive integer $N=N(d,h)$ such that the following holds.
Let $H$ be a graph on $h$ vertices of maximum degree $d$.
Let $G$ be a 4-edge-connected graph with no $H$-immersion.
Then $G$ is $(\chi_*(d-1)+1)$-colorable with clustering $N$.
Furthermore, if there exists exactly one vertex of $H$ having degree $d$, then $G$ is $(\chi_*(d-2)+1)$-colorable with clustering $N$.
\end{lemma}

\begin{pf}
Let $d \geq 3,h$ be positive integers.
We define the following.
	\begin{itemize}
		\item For every $x \in {\mathbb N}$, let $N_x$ be a positive integer such that every graph of maximum degree at most $x$ is $\chi_*(x)$-colorable with clustering $N_x$. 
		\item Let $\eta = \eta_{\ref{global}}(d,h)$ and $\xi = \xi_{\ref{global}}(d,h)$, where $\eta_{\ref{global}}$ and $\xi_{\ref{global}}$ are the integers $\eta$ and $\xi$ mentioned in Theorem \ref{global}.
		\item Define $N=N_{\ref{coloring_const_bag}}((N_{d-1}+N_{d-2})(\xi+1),\eta,h)$, where $N_{\ref{coloring_const_bag}}$ is the integer $N^*$ mentioned in Lemma \ref{coloring_const_bag}.
	\end{itemize}

Let $H$ be a graph on $h$ vertices of maximum degree $d$.
If $H$ has exactly one vertex of degree $d$, then let $\chi = \chi_*(d-2)$; otherwise, let $\chi = \chi_*(d-1)$.
Let $G$ be a 4-edge-connected graph with no $H$-immersion.
We shall prove that $G$ is $(\chi+1)$-colorable with clustering $N$.
Suppose to the contrary that $G$ is not $(\chi+1)$-colorable with clustering $N$.

Since $N \geq h$, $\lvert V(G) \rvert > \lvert V(H) \rvert$, for otherwise $G$ is 1-colorable with clustering $N$.
If $H$ is non-exceptional, then let $H'=H$; if $H$ is exceptional, then let $H'$ be a graph obtained from $H$ by subdividing an edge. 
Note that the maximum degree of $H'$ is $d$.
Since $d \geq 3$, if $H$ has exactly one vertex of degree $d$, then $H'$ has exactly one vertex of degree $d$.

Since $G$ does not contain an $H$-immersion, by Theorem \ref{global}, there exists a tree-cut decomposition $(T,\X)$ of $G$ of adhesion at most $\eta$ such that for every $t \in V(T)$, there exists $Z_t \subseteq E(G)$ with $\lvert Z_t \rvert \leq \xi$ such that if $G_t$ is the torso at $t$, then there exists a nonnegative integer $d_t \leq d$ such that 
	\begin{itemize}
		\item[(i)] the number of vertices of degree at least $d_t$ in $G_t-Z_t$ is less than the number of vertices of degree at least $d_t$ in $H'$, 
		\item[(ii)] every vertex of $G_t-Z_t$ of degree at least $d_t$ in $G_t-Z_t$ is a non-peripheral vertex of $G_t$,
		\item[(iii)] if $\lvert V(T) \rvert=1$ or $t$ is not a leaf of $T$, then every vertex in $X_t$ has degree at least $d_t$ in $G_t-Z_t$, and
		\item[(iv)] if $t$ is a leaf and $\lvert V(T) \rvert \geq 2$, then $\lvert X_t \rvert \leq 1$.
	\end{itemize}

If $\lvert V(T) \rvert=1$, then by (i) and (iii), $G$ contains at most $\lvert V(H') \rvert-1 \leq h$ vertices, so $G$ is 1-colorable with clustering $h \leq N$, a contradiction.
So $\lvert V(T) \rvert \geq 2$.
By (iv), $\lvert X_t \rvert \leq 1$ for every leaf $t$ of $T$.

For every $t \in V(T)$, since $d_t \leq d$, $(G_t-Z_t)-X_t$ has maximum degree at most $d_t-1 \leq d-1$ by (ii), so $(G_t-Z_t)-X_t$ is $\chi_*(d-1)$-colorable with clustering $N_{d-1}$.
By Lemma \ref{coloring_adding_edge}, for every $t \in V(T)$, $G_t-X_t$ is $\chi_*(d-1)$-colorable with clustering $(\lvert Z_t \rvert+1)N_{d-1} \leq (\xi+1) N_{d-1}$.
In addition, (i), (iii) and (iv) imply that $\lvert X_t \rvert \leq \max\{\lvert V(H') \rvert-1,1\} \leq h$ for every node $t \in V(T)$.
Since $d \geq 3$, $\chi^*(d-1) \geq 2$.
Hence by Lemma \ref{coloring_const_bag}, $G$ is $(\chi_*(d-1)+1)$-colorable with clustering $N_{\ref{coloring_const_bag}}((\xi+1) N_{d-1}, \eta, h) \leq N$.
	
So $\chi = \chi_*(d-2)$.
Hence $H$ has exactly one vertex of degree $d$.
So $H'$ has exactly one vertex of degree $d$.
Hence for every non-leaf $t$ of $T$, if $d_t = d$, then by (i) and (iii), $X_t=\emptyset$ and $G_t-Z_t$ has maximum degree at most $d-1$; if $d_t<d$, then the maximum degree of $(G_t-Z_t)-X_t$ is at most $d_t-1 \leq d-2$.
For every non-leaf $t$ of $T$, if $d_t=d$, then let $k_t = \chi_*(d-1)$; if $d_t<d$, then let $k_t = \max\{\chi_*(d-2),2-\lvert X_t \rvert\}$.
Then for every non-leaf $t$ of $T$, $(G_t-Z_t)-X_t$ is $k_t$-colorable with clustering $N_{d-1}+N_{d-2}$.
By Lemma \ref{coloring_adding_edge}, for every non-leaf $t$ of $T$, $G_t-X_t$ is $k_t$-colorable with clustering $(\xi+1) \cdot (N_{d-1}+N_{d-2})$.
Note that for every non-leaf $t$ of $T$, $k_t + \min\{\lvert X_t \rvert,1\}$ is at most either $\chi_*(d-1)$ or $\chi_*(d-2)+1$.
Since every graph with maximum degree at most $d-1$ can be partitioned into a stable set and a induced subgraph of maximum degree at most $d-2$, we have $\chi_*(d-1) \leq \chi_*(d-2)+1$.
So $\max_t\{k_t+\min\{\lvert X_t \rvert,1\}\} \leq \chi_*(d-2)+1$, where the maximum is over all non-leaves $t$ of $T$.
For every leaf $t$ of $T$, let $k_t=2-\lvert X_t \rvert$.
For every leaf $t$ of $T$, $\lvert X_t \rvert \leq 1$ by (iv), so $G_t$ is $k_t$-colorable with clustering $2 \leq N$.
Hence $\max_{t \in V(T)}\{k_t+\min\{\lvert X_t \rvert,1\}\} \leq \chi_*(d-2)+1$.
By Lemma \ref{coloring_const_bag}, $G$ is $(\chi_*(d-2)+1)$-colorable with clustering $N_{\ref{coloring_const_bag}}((N_{d-1}+N_{d-2})(\xi+1), \eta,h) = N$.
This proves the lemma.
\end{pf}

\bigskip

The following lemma is a simple variant of a result of Dirac \cite{d_cut}.
For every graph $G$ and subset $S$ of $V(G)$, we define $N_G(S) = \{v \in V(G)-S: v$ is adjacent to some vertex in $S\}$.

\begin{lemma} \label{coloring_cut}
Let $G$ be a graph.
Let $k$ and $N$ be positive integers.
Let $[A,B]$ be an edge-cut of $G$ of order at most $k-1$.
If both $G[A]$ and $G[B]$ are $k$-colorable with clustering $N$, then $G$ is $k$-colorable with clustering $N$.
\end{lemma}

\begin{pf}
Let $c_A$ and $c_B$ be $k$-colorings of $G[A]$ and $G[B]$ with clustering $N$, respectively.
Define $H$ to be a simple bipartite graph $H$ with $V(H)=\{a_i,b_i: i \in [k]\}$ and with bipartition $(\{a_i: i \in [k]\},\{b_i: i \in [k]\})$ such that two vertices $a_i$ and $b_j$ are adjacent in $H$ if and only if there exists an edge incident with a vertex $v \in A$ with $c_A(v)=i$ and a vertex $u \in B$ with $c_B(u)=j$.
Note that $\lvert E(H) \rvert \leq \lvert [A,B] \rvert \leq k-1$.
Let $H'$ be the bipartite complement of $H$.
That is, $V(H')=V(H)$ and $E(H') = \{a_ib_j: i,j \in [k]\}-E(H)$.

Suppose that $H'$ does not contain a perfect matching.
Then by Hall's theorem, there exists $S \subseteq \{a_i: i \in [k]\}$ such that $\lvert N_{H'}(S) \rvert < \lvert S \rvert$.
Note that every vertex in $S$ is adjacent in $H$ to every vertex in $\{b_j: j \in [k]\}-N_{H'}(S)$.
Hence $H$ contains at least $\lvert S \rvert \cdot (k-\lvert N_{H'}(S) \rvert) \geq \lvert S \rvert (k-\lvert S \rvert+1)$ edges.
Since $\lvert N_{H'}(S) \rvert < \lvert S \rvert$, $S \neq \emptyset$. 
Hence $1 \leq \lvert S \rvert \leq k$.
So $H$ contains at least $\lvert S \rvert (k-\lvert S \rvert+1) \geq k$ edges, a contradiction.

Hence $H'$ has a perfect matching $\{a_ib_{\sigma(i)}: i \in [k]\}$ for some bijection $\sigma: [k] \rightarrow [k]$.
So for each $i \in [k]$, $a_i$ is not adjacent to $b_{\sigma(i)}$ in $H$.
Let $c_A'$ be the $k$-coloring of $G[A]$ such that for every $i \in [k]$ and $v \in A$ with $c_A(v)=i$, $c_A'(v)=\sigma(i)$.
Note that $c'_A$ is a $k$-coloring of $G[A]$ with clustering $N$.
Define $c$ to be the $k$-coloring such that for every $v \in A$, $c(v)=c'_A(v)$, and for every $v \in B$, $c(v)=c_B(v)$.
Then each $c$-monochromatic component of $G$ is contained in $G[A]$ or $G[B]$.
So $c$ is a $k$-coloring of $G$ with clustering $N$.
\end{pf}

\bigskip

Recall that for every graph $H$, $\chi_*(H)$ is the minimum $k$ such that there exists $N \in {\mathbb N}$ such that every graph with no $H$-immersion is $k$-colorable with clustering $N$.

\begin{lemma} \label{d=1}
If $H$ is a graph of maximum degree 1, then $\chi_*(H)=1$.
\end{lemma}

\begin{pf}
Let $N=(\lvert V(H) \rvert-1)^{\lvert V(H) \rvert}$.

Let $G$ be a graph with no $H$-immersion.
Since $H$ has maximum degree 1, $H$ is a disjoint union of copies of $K_2$ and isolated vertices.
Note that $K_{1,\lvert V(H) \rvert}$ and the path on $\lvert V(H) \rvert$ vertices contain an $H$-immersion.
Hence the maximum degree of $G$ is at most $\lvert V(H) \rvert-1$, and every path in $G$ contains at most $\lvert V(H) \rvert-1$ vertices.
So every component of $G$ contains at most $(\lvert V(H) \rvert-1)^{\lvert V(H) \rvert}=N$ vertices.
Therefore, $G$ is 1-colorable with clustering $N$.
\end{pf}

\bigskip

Now we are ready to prove Theorem \ref{main_intro}.
The following is a restatement of Theorem \ref{main_intro}.

\begin{theorem} \label{main_clu}
Let $d$ be a positive integer, and let $H$ be a graph of maximum degree $d$.
	\begin{enumerate}
		\item If $d=1$, then $\chi_*(H)=1$.
		\item If $d \geq 2$ and $H$ has exactly one vertex of degree $d$, then $\chi_*(d-1) \leq \chi_*(H) \leq \max\{\chi_*(d-2)+1,4\}$.
		\item If $d \geq 2$ and $H$ has at least two vertices of degree $d$, then $\chi_*(d-2)+1 \leq \chi_*(H) \leq \max\{\chi_*(d-1)+1,4\}$.
	\end{enumerate}
\end{theorem}

\begin{pf}
Statement 1 immediate follows from Lemma \ref{d=1}.
So we may assume $d \geq 2$.

Let $\chi = \max\{\chi_*(d-2)+1,4\}$ if $H$ has exactly one vertex of degree $d$; otherwise, let $\chi = \max\{\chi_*(d-1)+1,4\}$.
Let $N=N_{\ref{coloring_upper_bdd}}(d,\lvert V(H) \rvert)$, where $N_{\ref{coloring_upper_bdd}}$ is the number $N$ mentioned in Lemma \ref{coloring_upper_bdd}.

We first prove the upper bounds.

Suppose that $d \geq 3$ and there exists a graph $G$ with no $H$-immersion such that $G$ is not $\chi$-colorable with clustering $N$.
We further assume that $\lvert V(G) \rvert$ is as small as possible.
By Lemma \ref{coloring_upper_bdd}, $G$ is not 4-edge-connected.
So there exists an edge-cut $[A,B]$ of $G$ of order at most 3 with $A \neq \emptyset \neq B$.
Note that both $G[A]$ and $G[B]$ are subgraphs of $G$, so they do not contain an $H$-immersion.
By the minimality of $G$, $G[A]$ and $G[B]$ are $\chi$-colorable with clustering $N$.
Since $\chi \geq 4 > \lvert [A,B] \rvert$, by Lemma \ref{coloring_cut}, $G$ is $\chi$-colorable with clustering $N$, a contradiction.

Therefore, if $d \geq 3$, then every graph with no $H$-immersion is $\chi$-colorable with clustering $N$.
If $d \leq 2$, then let $H'$ be the graph obtained from $H$ by adding a loop incident with a vertex of degree $d$, so $H'$ has maximum degree $d+2$ with $3 \leq d+2 \leq 4$, and hence $\chi_*(H) \leq \chi_*(H') \leq \max\{\chi_*(3)+1,4\}=4 \leq \chi$.
This proves the upper bound for Statements 2 and 3 of this theorem.

Now we prove the lower bounds.

Every graph of maximum degree at most $d-1$ does not contain an $H$-immersion.
So $\chi_*(H) \geq \chi_*(d-1)$.
This proves Statement 2.
To prove Statement 3, it suffices to show that $\chi_*(d-2)+1 \leq \chi_*(H)$ when $H$ has at least two vertices of degree $d$.

Now we assume that $H$ has at least two vertices of degree $d$.
Suppose to the contrary that $\chi_*(H) \leq \chi_*(d-2)$.
So there exists a positive integer $\eta$ such that every graph with no $H$-immersion is $\chi_*(d-2)$-colorable with clustering $\eta$.
By the definition of $\chi_*(d-2)$, there exists a graph $L$ of maximum degree at most $d-2$ such that there exists no $(\chi_*(d-2)-1)$-coloring of $L$ with clustering $\eta$.

Define $Q$ to be the graph obtained from a union of $\eta$ disjoint copies $L_1,L_2,...,L_\eta$ of $L$ by adding a vertex $v^*$ adjacent to all other vertices.
Since $L$ is of maximum degree at most $d-2$, $Q$ has at most one vertex of degree at least $d$.
Since $H$ contains at least two vertices of degree $d$, $Q$ does not contain an $H$-immersion.
So there exists a $\chi_*(d-2)$-coloring $c$ of $Q$ with clustering $\eta$.
By symmetry, we may assume that $c(v^*)=\chi_*(d-2)$.
Since $c$ is of clustering $\eta$, there exists $i \in [\eta]$ such that $c(v) \neq c(v^*)$ for every $v \in V(L_i)$.
So the restriction of $c$ on $L_i$ is a $(\chi_*(d-2)-1)$-coloring of clustering $\eta$.
However, it is impossible by the definition of $L$.
This proves the theorem.
\end{pf}

\section{Application to tree-decompositions} \label{sec:tree_decomp}

Let $G$ be a graph.
A {\it tree-decomposition} of $G$ is a pair $(T,\X)$ such that $T$ is a tree and $\X$ is a collection $(X_t: t \in V(T))$ of subsets of $V(G)$ such that
	\begin{itemize}
		\item $\bigcup_{t \in V(T)}X_t = V(G)$,
		\item for every $e \in E(G)$, there exists $t \in V(T)$ such that $X_t$ contains all ends of $e$, and 
		\item for every $v \in V(G)$, the set $\{t \in V(T): v \in X_t\}$ induces a connected subgraph of $T$.
	\end{itemize}
The {\it adhesion} of $(T,\X)$ is $\max_{tt' \in E(T)}\lvert X_t \cap X_{t'} \rvert$.
For every $t \in V(T)$, the {\it torso} at $t$ in $(T,\X)$ is the graph obtained from $G[X_t]$ by for each neighbor $t'$ of $t$, adding edges such that $X_t \cap X_{t'}$ is a clique. 

\begin{lemma} \label{simply_td}
Let $d$ and $\eta$ be positive integers.
Let $G$ be a graph with maximum degree at most $d$.
Let $(T,\X)$ be a tree-decomposition of $G$ of adhesion at most $\eta$.
Then there exists a tree-decomposition $(T,\X'=(X'_t: t \in V(T)))$ of $G$ such that 
	\begin{enumerate}
		\item the adhesion of $(T,\X')$ is at most the adhesion of $(T,\X)$.
		\item for every $t \in V(T)$, $X'_t \subseteq X_t$, and the torso at $t$ in $(T,\X')$ is a subgraph of the torso at $t$ in $(T,\X)$,
		\item for every $t \in V(T)$ and $v \in X'_t$, there exist at most $d+1$ neighbors $t'$ of $t$ such that $v \in X'_{t'}$, and
		\item for every $t \in V(T)$, the torso at $t$ in $(T,\X')$ has maximum degree at most $\eta d+\eta-1$.
	\end{enumerate}
\end{lemma}

\begin{pf}
Let $r$ be a node of $T$.
We assume that $T$ is a rooted tree rooted at $r$.
For each $t \in V(T)$, let $T_t$ be the maximal subtree of $T$ rooted at $t$.
For every $t \in V(T)$ and every vertex $v \in X_t$, define $C_{t,v}$ to be the set of children $c$ of $t$ such that $(\bigcup_{t' \in V(T_c)}X_{t'})-X_t$ contains a neighbor of $v$.
Since the maximum degree of $G$ is at most $d$, $\lvert C_{t,v} \rvert \leq d$. 
For every $v \in V(G)$, let $r_v$ be the node of $T$ with $v \in X_{r_v}$ closest to $r$.

Denote $\X$ by $(X_t: t \in V(T))$.
For every $t \in V(T)$, define $X'_t = \{v \in X_t:$ either $t = r_v$, or $t \in C_{p,v}$, where $p$ is the parent of $t\}$.
Then $(T,\X')$ is a tree-decomposition satisfying Statements 1-3.
By Statement 3, for every $t \in V(T)$, the torso at $t$ in $(T,\X')$ has maximum degree at most $d+(\eta-1)(d+1) = \eta d+\eta-1$.
So Statement 4 holds.
\end{pf}

\begin{lemma} \label{tree-cut_tree}
Let $d$ and $\eta$ be positive integers.
Let $G$ be a graph with maximum degree at most $d$.
Let $(T,\X)$ be a tree-decomposition of $G$ of adhesion at most $\eta$.
Denote $\X$ by $(X_t: t \in V(T))$.
Then there exists a tree-cut decomposition $(T',\X' = (X'_t: t \in V(T')))$ of $G$ of adhesion at most $(d+1)^2\eta+d$ such that the following hold.
	\begin{enumerate}
		\item $T'$ is obtained from $T$ by attaching leaves. 
		\item For every $t \in V(T')$, $\lvert X'_t \rvert \leq 1$, and if $\lvert X'_t \rvert = 1$, then $t$ is a leaf of $T'$.
		\item For every $t \in V(T)$, the torso at $t$ in $(T',\X')$ is a subgraph of a graph $R_t$ obtained from a subgraph of the torso at $t$ in $(T,\X)$ by identifying a set of at most $\eta$ vertices in $X_t$ into a vertex, adding a set $I$ of vertices and adding edges incident with $I$ such that $I$ is a stable set and the neighborhood of every vertex in $I$ is a clique of size at most $\eta$.
		\item For every $t \in V(T)$, the maximum degree of $R_t$ is at most $(d+1)\eta^2+d$. 
	\end{enumerate}
\end{lemma}

\begin{pf}
By Lemma \ref{simply_td}, there exists a tree-decomposition $(T,\X^1=(X^1_t: t \in V(T)))$ of $G$ of adhesion at most $\eta$ such that for every $t \in V(T)$, the torso at $t$ in $(T,\X^1)$ is a subgraph of the torso at $t$ in $(T,\X)$, and for every $v \in X^1_t$, there exist at most $d+1$ neighbors $t'$ of $t$ such that $v \in X^1_{t'}$.
For every $t \in V(T)$, let $Q_t$ be the torso at $t$ in $(T,\X^1)$.

Let $r$ be a node of $T$.
We assume that $T$ is a rooted tree rooted at $r$.
For every $t \in V(T)-\{r\}$, let $p_t$ be the parent of $t$.
For every vertex $v$ of $G$, let $r_v$ be the node of $T$ with $v \in X^1_{r_v}$ closest to $r$.
For every $t \in V(T)$, define $X^2_t = \{v \in X^1_t: t=r_v\}$.
Let $\X^2 = (X^2_t: t \in V(T))$.
Then $(T,\X^2)$ is a tree-cut decomposition of $G$.

For every $tt' \in E(T)$, if $uv$ is an edge of $G$ contained in $\adh_{(T,\X^2)}(tt')$, then $\{u,v\} \cap X^1_t \cap X^1_{t'} \neq \emptyset$.
Since the maximum degree of $G$ is at most $d$, the adhesion of $(T,\X^2)$ is at most $d\eta$. 

Note that for every $t \in V(T)$, if there exists $v \in X^1_t-X^2_t$, then $t \neq r$ and $v \in X^1_t \cap X^1_{p_t}$.
Let $S_r = \emptyset$, and for every $t \in V(T)-\{r\}$, let $S_t = X^1_t \cap X^1_{p_t}$.
So for every $t \in V(T)-\{r\}$, the peripheral vertex of the torso at $t$ in $(T,\X^2)$ corresponding to the component of $T-t$ containing $r$ is obtained from $Q_t$ by identifying $S_t$ into a vertex $v_t$ and deleting the resulting loops.
Note that $\lvert S_t \rvert \leq \eta$, since the adhesion of $(T,\X^1)$ is at most $\eta$.

Suppose that there exist $t \in V(T)$, a non-loop edge $uv$ of $G$, and two distinct components $T_1,T_2$ of $T-t$ not containing $r$ such that $u \in X^2_{T_1}$ and $v \in X^2_{T_2}$. 
Then $r_u \in V(T_1)$ and $r_v \in V(T_2)$.
So there exists no $t' \in V(T)$ such that $X^1_{t'} \supseteq \{u,v\}$, contradicting that $(T,\X^1)$ is a tree-decomposition.

Hence for every $t \in V(T)$, the peripheral vertices of the torso at $t$ in $(T,\X^2)$ corresponding to the components of $T-t$ disjoint from $r$ form a stable set, denoted by $I_t$.
In addition, since $(T,\X^1)$ is a tree-decomposition of $G$, for every $t \in V(T)$, if $q \in I_t$, then the neighborhood of $q$ in the torso at $t$ in $(T,\X^2)$ is contained in $(X^1_t \cap \bigcup_{t' \in V(W_q)}X^1_{t'} - S_t) \cup \{v_t\}$, where $W_q$ is the component of $T-t$ corresponding to $q$; and if $q$ is adjacent to $v_t$, then $S_t \cap X^1_t \cap \bigcup_{t' \in V(W_q)}X^1_{t'} \neq \emptyset$.

For every $t \in V(T)$, let $R^2_t$ be the graph obtained from $Q_t$ by identifying $S_t$ into a vertex $v_t$ and adding $I_t$ and edges such that for every $q \in I_t$, the neighborhood of $q$ in $R^2_t$ is the same as the neighborhood of $q$ in the torso at $t$ in $(T,\X^2)$.
Since the adhesion of $(T,\X^1)$ is at most $\eta$, the neighborhood in $R^2_t$ of each vertex in $I_t$ is a clique of size at most $\eta$.
Note that the torso at $t$ in $(T,\X^2)$ is a subgraph of $R^2_t$. 

Since for every $t \in V(T)$ and $v \in X^1_t$, there exist at most $d+1$ neighbors $t'$ of $t$ such that $v \in X^1_{t'}$, the maximum degree of $R^2_t$ is at most $(d+(d+1)\eta)\eta \leq (d+1)\eta^2+d$.

Define $T'$ to be the tree obtained from $T$ by for each $t \in V(T)$, attaching $\lvert X^2_t \rvert$ leaves adjacent to $t$.
So for every $t \in V(T)$, there exists a bijection $\sigma_t$ from the set $L_t$ of leaves attached on $t$ to $X^2_t$.
For each $t \in V(T)$, define $X'_t = \emptyset$; for each $t \in V(T')-V(T)$, there uniquely exists $t' \in V(T)$ such that $t \in L_{t'}$, and we define $X'_t = \{\sigma_{t'}(t)\}$.
Define $\X' = (X'_t: t \in V(T'))$.
Then $(T',\X')$ is a tree-cut decomposition of $G$ such that Statements 1 and 2 of this lemma hold.

For every $t \in V(T)$, define $R_t$ to be a graph obtained from $R^2_t$ by deleting all loops.
For every $t \in V(T)$, since the torso at $t$ in $(T',\X')$ is obtained from the torso at $t$ in $(T,\X^2)$ by deleting all loops, the torso at $t$ in $(T',\X')$ is a subgraph of $R_t$.
Furthermore, the maximum degree of $R_t$ is at most the maximum degree of $R^2_t$.
Hence Statements 3 and 4 hold.

Since the adhesion of $(T',\X')$ is at most the maximum degree of the torsos in $(T',\X')$, the adhesion of $(T',\X')$ is at most $(d+1)\eta^2+d$.
This proves the lemma.
\end{pf}

\begin{lemma} \label{color_id_add}
For any positive integers $d, \eta, N,d'$, there exists an integer $N^*=N^*(d,\eta,N,d')$ such that the following hold.
Let $k$ be a positive integer, and let $G$ be a graph with maximum degree at most $d$ such that $G$ is $k$-colorable with clustering $N$.
Let $G'$ be a graph with maximum degree at most $d'$ obtained from $G$ by identifying a set of at most $\eta$ vertices into a vertex, adding a set $I$ of vertices and adding edges incident with $I$ such that $I$ is a stable set in $G'$, and the neighborhood of each vertex in $I$ is a clique of size at most $\eta$.
Then $G'$ is $k$-colorable with clustering $N^*$.
\end{lemma}

\begin{pf}
Let $d,\eta,N,d'$ be positive integers.
Let $N_0 = d\eta N+1$.
Define $N^* = (d'+1)N_0$.

Let $k$ be a positive integer.
Let $G$ be a graph with maximum degree at most $d$ such that $G$ is $k$-colorable with clustering $N$.
Let $G_0$ be a graph with maximum degree at most $d'$ obtained from $G$ by identifying a set $S$ of at most $\eta$ vertices into a vertex $v_S$.
Let $G'$ be a graph with maximum degree at most $d'$ obtained from $G_0$ by adding a set $I$ of vertices and adding edges incident with $I$ such that $I$ is a stable set in $G'$, and the neighborhood of each vertex in $I$ is a clique of size at most $\eta$.
It suffices to prove that $G'$ is $k$-colorable with clustering $N^*$.

Let $f$ be a $k$-coloring of $G$ with clustering $N$.
Let $f_0(v_S)=1$.
For every $v \in V(G_0)-\{v_S\}$, let $f_0(v)=f(v)$.
Since the maximum degree of $G$ is at most $d$, $f_0$ is a $k$-coloring of $G_0$ with clustering $d\lvert S \rvert \cdot N+1 \leq d \eta N+1 = N_0$.

For every $v \in V(G_0)$, let $f'(v)=f_0(v)$; for every $v \in I$, let $f'(v)=1$.
Then $f'$ is a $k$-coloring of $G'$.
Let $M$ be an $f'$-monochromatic component of $G'$.
It suffices to show that $\lvert V(M) \rvert \leq N^*$.

If $V(M) \subseteq I$, then since $I$ is a stable set in $G'$, $\lvert V(M) \rvert =1$.
So we may assume that $V(M-I) \neq \emptyset$.
Since $I$ is a stable set in $G'$, and the neighborhood of each vertex in $I$ is a clique in $G'$, $M-I$ is connected.
So $M-I$ is a $f_0$-monochromatic component of $G_0$.
Hence $\lvert V(M-I) \rvert \leq N_0$.
Since the maximum degree of $G'$ is at most $d'$, and $I$ is a stable set in $G'$, $\lvert V(M) \cap I \rvert \leq d' \lvert V(M-I) \rvert \leq d'N_0 $.
So $\lvert V(M) \rvert \leq \lvert V(M) \cap I \rvert + \lvert V(M-I) \rvert \leq (d'+1)N_0$.
This proves the lemma.
\end{pf}

\bigskip

The following is a restatement of Statement 1 in Corollary \ref{td_intro}.

\begin{theorem}
For any positive integers $\eta,d$ and $N$, there exists a positive integer $N^*$ such that for every positive integer $k \geq 2$, if $G$ is a graph with maximum degree at most $d$ and $G$ admits a tree-decomposition $(T,\X)$ of adhesion at most $\eta$ such that for every $t \in V(T)$, the torso at $t$ in $(T,\X)$ is $k$-colorable with clustering $N$, then $G$ is $k$-colorable with clustering $N^*$.
\end{theorem}

\begin{pf}
Let $\eta,d,N$ be positive integers.
Let $N_1 = N_{\ref{color_id_add}}(\eta d+\eta-1,\eta,N,(d+1)\eta^2+d)$, where $N_{\ref{color_id_add}}$ is the integer $N^*$ mentioned in Lemma \ref{color_id_add}.
Define $N^*=N_{\ref{coloring_const_bag}}(N_1,(d+1)^2\eta+d,1)$, where $N_{\ref{coloring_const_bag}}$ is the integer $N^*$ mentioned in Lemma \ref{coloring_const_bag}.

Let $k$ be an integer with $k \geq 2$.
Let $G$ be a graph with maximum degree at most $d$ such that $G$ admits a tree-decomposition $(T,\X)$ of adhesion at most $\eta$ such that for every $t \in V(T)$, the torso at $t$ in $(T,\X)$ is $k$-colorable with clustering $N$.

Let $\C_0$ be the collection consisting of the graphs that are subgraphs of a torso at $t$ in $(T,\X)$ for some $t \in V(T)$.
By assumption, every graph in $\C_0$ is $k$-colorable with clustering $N$.
By Lemma \ref{simply_td}, there exists a tree-decomposition $(T,\X^0=(X^0_t: t \in V(T)))$ of $G$ of adhesion at most $\eta$ such that for every $t \in V(T)$, the torso at $t$ in $(T,\X^0)$ belongs to $\C_0$ and has maximum degree at most $\eta d+\eta-1$. 

Let $\C_0'$ be the set of graphs in $\C_0$ with maximum degree at most $\eta d+\eta-1$.
Let $\C_1$ be the collection consisting of the graphs of maximum degree at most $(d+1)\eta^2+d$ that can be obtained from some graph in $\C'_0$ by identifying a set of at most $\eta$ vertices into a vertex, adding a set $I$ of vertices and adding edges incident with $I$ such that $I$ is a stable set and the neighborhood of each vertex in $I$ is a clique of size at most $\eta$.
By Lemma \ref{color_id_add}, every graph in $\C_1$ is $k$-colorable with clustering $N_1$,

By Lemma \ref{tree-cut_tree}, there exists a tree-cut decomposition $(T',\X'=(X'_t: t\in V(T')))$ of $G$ of adhesion at most $(d+1)\eta^2+d$ such that the following hold.
	\begin{itemize}
		\item $T'$ is obtained from $T$ by adding leaves.
		\item For every $t \in V(T')$, $\lvert X'_{t'} \rvert \leq 1$, and if $\lvert X'_{t'} \rvert = 1$, then $t$ is a leaf of $T'$.
		\item For every $t \in V(T)$, the torso at $t$ in $(T',\X')$ has maximum degree at most $(d+1)\eta^2+d$ and is a subgraph of a graph in $\C_1$. 
		\item For every $t \in V(T')-V(T)$, the torso at $t$ in $(T',\X')$ has at most 2 vertices.
	\end{itemize}
For every $t \in V(T')$, if $X'_t \neq \emptyset$, then let $k_t=1$; otherwise let $k_t=k$.
So if $X'_t \neq \emptyset$, then $t$ is a leaf of $T'$, $k_t+\lvert X'_t \rvert=2$, and the torso at $t$ in $(T',\X')$ is $(k_t+\lvert X'_t \rvert)$-colorable with clustering $1 \leq N_1$; if $X'_t = \emptyset$, then since the torso at $t$ in $(T',\X')$ either is in $\C_1$ or has at most two vertices, it is $k_t$-colorable with clustering $N_1$.
Hence by Lemma \ref{coloring_const_bag}, $G$ is $\max_{t \in V(T')}\{k_t+\min\{\lvert X'_t \rvert,1\}\}$-colorable with clustering $N^*$.
Note that $\max_{t \in V(T')}\{k_t+\min\{\lvert X'_t \rvert,1\}\} \leq \max\{k,2\} = k$.
This proves the theorem.
\end{pf}

\bigskip

A similar and simpler argument proves the following theorem which is a restatement of Statement 2 in Corollary \ref{td_intro}.

\begin{theorem}
For any positive integers $\eta,d$ and $N$, there exists a positive integer $N^*$ such that for every integer $k$, if $G$ is a graph with maximum degree at most $d$ and $G$ admits a tree-decomposition $(T,\X=(X_t: t \in V(G)))$ of adhesion at most $\eta$ such that for every $t \in V(T)$, $G[X_t]$ is $k$-colorable with clustering $N$, then $G$ is $(k+1)$-colorable with clustering $N^*$.
\end{theorem}

\begin{pf}
Let $\eta,d,N$ be positive integers.
Let $N_1 = N_{\ref{color_id_add}}(\eta d+\eta-1,\eta,N,(d+1)\eta^2+d)$, where $N_{\ref{color_id_add}}$ is the integer $N^*$ mentioned in Lemma \ref{color_id_add}.
Define $N^*=N_{\ref{coloring_const_bag}}(N_1,(d+1)^2\eta+d,1)$, where $N_{\ref{coloring_const_bag}}$ is the integer $N^*$ mentioned in Lemma \ref{coloring_const_bag}.

Let $k$ be an integer and $G$ a graph as stated in this theorem.
Since $G[X_t]$ is $k$-colorable for every $t \in V(T)$, $k \geq 1$.

By Lemma \ref{simply_td}, there exist a set $\C_0$ of graphs that are $k$-colorable with clustering $N$ and a tree-decomposition $(T,\X^0=(X^0_t: t \in V(T)))$ of $G$ of adhesion at most $\eta$ such that for every $t \in V(T)$, $G[X^0_t] \in \C_0$ and has maximum degree at most $\eta d+\eta-1$. 
Let $\C_1$ be the collection consisting of the graphs of maximum degree at most $(d+1)\eta^2+d$ that can be obtained from some graph in $\C_0$ by identifying a set of at most $\eta$ vertices into a vertex. 
By Lemma \ref{color_id_add}, every graph in $\C_1$ is $k$-colorable with clustering $N_1$,

Let $\C_2$ be the collection consisting of the graphs of maximum degree at most $(d+1)^2\eta+d$ that can be obtained from some graph in $\C_1$ by adding a set $I$ of vertices and edges incident with $I$ such that $I$ is a stable set.
Then every graph in $\C_2$ is $(k+1)$-colorable with clustering $N_1$ since we can use a new color to color $I$.

By Lemma \ref{tree-cut_tree}, there exists a tree-cut decomposition $(T',\X'=(X'_t: t \in V(T')))$ of $G$ of adhesion at most $(d+1)\eta^2+d$ such that 
	\begin{itemize}
		\item for every $t \in V(T')$, $\lvert X'_t \rvert \leq 1$, and if $\lvert X'_t \rvert=1$, then $t$ is a leaf in $T'$, 
		\item for every $t \in V(T)$, the torso at $t$ in $(T',\X')$ is in $\C_2$, and
		\item for every $t \in V(T')-V(T)$, the torso at $t$ has at most 2 vertices.
	\end{itemize}
For every $t \in V(T')$, let $k_t=(1-\lvert X'_t \rvert)k+1$, so $k_t+\lvert X'_t \rvert \geq 2$.
For $t \in V(T')$, if $\lvert X'_t \rvert=0$, then $k_t=k+1$, so the torso at $t$ in $(T',\X')$ is $k_t$-colorable with clustering $N_1$; if $\lvert X'_t \rvert>0$, then $t$ is a leaf in $T'$, so the torso at $t$ in $(T',\X')$ has at most 2 vertices and is $k_t$-colorable with clustering $2 \leq N_1$.
By Lemma \ref{coloring_const_bag}, $G$ is $\max_{t \in V(T')}\{k_t+\min\{\lvert X'_t \rvert,1\}\}$-colorable with clustering $N^*$. 
Since $\max_{t \in V(T')}\{k_t+\min\{\lvert X'_t \rvert,1\}\} \leq k+1$, this proves the theorem.
\end{pf}

\section{Concluding remarks} \label{sec:remarks}

In this paper we prove that $\chi_*(H)$ is very close to $\chi_*(\Delta(H)-1)$ for every graph $H$.
But it remains unclear what $\chi_*(x)$ is, even for its asymptotic behavior.
It can be shown that $\lim_{x \rightarrow \infty} \frac{\chi_*(x)}{x}$ exists by a result of Lov\'{a}sz \cite{l_partition}.
And $\frac{1}{4} \leq \lim_{x \rightarrow \infty} \frac{\chi_*(x)}{x} \leq \frac{1}{3}$ by \cite{adov,hst}.

\begin{question}
Determine $\lim_{x \rightarrow \infty}\frac{\chi_*(x)}{x}$.
\end{question}

Another natural question is about strong immersion.
The immersion containment can be equivalently defined as follows.
For graphs $G$ and $H$, we say that $G$ contains an {\it $H$-immersion} if there exist functions $\pi_V$ and $\pi_E$ such that 
	\begin{itemize}
		\item $\pi_V$ is an injection from $V(H)$ to $V(G)$,
		\item $\pi_E$ maps each edge of $H$ to a subgraph of $G$ such that for each $e \in E(H)$, if $e$ has distinct ends $x,y$, then $\pi_E(e)$ is a path in $G$ with ends $\pi_V(x)$ and $\pi_V(y)$, and if $e$ is a loop with end $v$, then $\pi_E(e)$ is a cycle containing $\pi_V(v)$, and
		\item if $e_1,e_2$ are distinct edges of $H$, then $\pi_E(e_1)$ and $\pi_E(e_2)$ are edge-disjoint.
	\end{itemize}
We say that a graph $G$ contains another graph $H$ as a {\it strong immersion} if $G$ contains an $H$-immersion such that the witness functions $\pi_V$ and $\pi_E$ satisfy the extra property that for every $v \in V(H)$ and $e \in E(H)$, if $e$ is not incident with $v$, then $\pi_E(e)$ does not contain $\pi_V(v)$.

Strong immersion was introduced by Nash-Williams, and numerous problems that were proposed for minors, topological minors and immersions have been proposed for strong immersion as well.
So it is natural to consider the clustered chromatic number of the class of graphs with no $H$-strong immersion for any fixed graph $H$.
It turns out that the answer is different from the one for immersion for some graph $H$, but possibly not too much.

The {\it clustered chromatic number} of a class $\C$ of graphs is the minimum $k$ such that there exists a positive integer $N$ such that every graph in $\C$ is $k$-colorable with clustering $N$.

\begin{proposition} \label{strong_imm}
Let $d$ be a positive integer.
Let $H$ be a graph with maximum degree $d$ such that there exists a cycle in $H$ containing at least 3 vertices of degree $d$.
Then the clustered chromatic number of the class of graphs that do not contain $H$ as a strong immersion is at least $\chi_*(d-3)+2$.
\end{proposition}

\begin{pf}
Suppose to the contrary that there exists a positive integer $N$ such that every graph that does not contain $H$ as a strong immersion is $(\chi_*(d-3)+1)$-colorable with clustering $N$.
Let $L$ be a graph with maximum degree at most $d-3$ such that no $(\chi_*(d-3)-1)$-coloring of $L$ with clustering $N$ exists.
Let $L^*$ be the simple graph obtained from a path $v_1v_2...v_{N+1}$ on $N+1$ vertices by for each $i \in [N]$, adding $2N-1$ disjoint copies $L_{i,1},L_{i,2},...,L_{i,2N-1}$ of $L$ and adding edges such that $v_i$ and $v_{i+1}$ are adjacent to all vertices in $\bigcup_{j=1}^{2N-1}L_{i,j}$.

Since $L$ has maximum degree at most $d-3$, $v_1,v_2,...,v_{N+1}$ are the only vertices in $L^*$ with degree at least $d$.
Suppose that $L^*$ contains $H$ as a strong immersion.
Since there exists a cycle in $H$ containing 3 vertices of degree $d$, there exist distinct elements $\alpha<\beta<\gamma$ in $[N+1]$ and edge-disjoint paths $P_1,P_2,P_3$ in $L^*$ such that each $P_i$ contains exactly 2 vertices in $\{v_\alpha,v_\beta,v_\gamma\}$.
But there exists no path in $L^*-v_\beta$ from $v_\alpha$ to $v_\gamma$, contradicting the existence of $P_1,P_2,P_3$.

Hence $L^*$ does not contain $H$ as a strong immersion.
So by assumption, there exists a $(\chi_*(d-3)+1)$-coloring $c$ of $L^*$ with clustering $N$.
So the path $v_1v_2,,,v_{N+1}$ on $N+1$ vertices is not $c$-monochromatic.
Hence there exists $i \in [N]$ such that $c(v_i) \neq c(v_{i+1})$.
By symmetry, we may assume that $c(v_i)=\chi_*(d-3)$ and $c(v_{i+1})=\chi_*(d-3)+1$.
Since $c$ has clustering $N$, there are most $N-1$ indices $j$ such that $L_{i,j}$ contains a vertex with color $\chi_*(d-3)$ and at most $N-1$ indices $j'$ such that $L_{i,j'}$ contains a vertex with color $\chi_*(d-3)+1$.
Hence there exists an index $j^* \in [2N-1]$ such that the restriction of $c$ on $L_{i,j^*}$ is a $(\chi_*(d-3)-1)$-coloring with clustering $N$, contradicting the definition of $L$.
\end{pf}

\bigskip

Note that there are infinitely many positive integers $d$ such that $\chi_*(d-3)+2 > \chi_*(d-1)+1$, for otherwise $\lim_{x \rightarrow \infty} \frac{\chi_*(x)}{x} \geq \frac{1}{2}$, contradicting $\lim_{x \rightarrow \infty} \frac{\chi_*(x)}{x} \leq \frac{1}{3}$.
Hence Proposition \ref{strong_imm} and Theorem \ref{main_clu} show that the clustered chromatic number of $H$-immersion free graphs and $H$-strong immersion free graphs are different for infinitely many graphs $H$.
However, it is unknown whether the gap can be arbitrarily large.
We conjecture that it is not the case.

\begin{conjecture}
There exists a positive integer $C$ such that for every graph $H$, the clustered chromatic number of the class of graphs that do not contain $H$ as a strong immersion is at most $\chi_*(H)+C$.
\end{conjecture}


\begin{thebibliography}{99}

\bibitem{al}
F. N. Abu-Khzam and M. A. Langston, {\it Graph coloring and the immersion order}, Computing and Combinatorics, Lecture Notes in Computer Science, Vol. 2697 (2003), 394--403.

\bibitem{adov}
N. Alon, G. Ding, B. Oporowski and D. Vertigan, {\it Partitioning into graphs with only small components}, J. Combin. Theory Ser. B 87 (2003), 231--243.

\bibitem{ah}
K. Appel and W. Haken, {\it Every planar map is four colorable. I. Discharging}, Illinois J. Math. 21 (1977), 429--490.

\bibitem{ahk}
K. Appel, W. Haken and J. Koch, {\it Every planar map is four colorable. II. Reducibility}, Illinois J. Math. 21 (1977), 491--567.

\bibitem{c_hajos}
P. A. Catlin, {\it Haj\'{o}s' graph-coloring conjecture: variations and counterexamples}, J. Combin. Theory Ser. B 26 (1979), 268--274.

\bibitem{dp}
M. Delcourt and L. Postle, {\it Reducing linear Hadwiger's conjecture to coloring small graphs}, arXiv:2108.01633.

\bibitem{ddfmms}
M. DeVos, Z. Dvo\v{r}\'{a}k, J. Fox, J. McDonald, B. Mohar and D. Scheide, {\it Minimum degree condition forcing complete graph immersion}, Combinatorica 34 (2014), 279--298.

\bibitem{dkmo}
M. DeVos, K. Kawarabayashi, B. Mohar and H. Okamura, {\it Immersing small complete graphs}, Ars Math. Contemp. 3 (2010), 139--146.

\bibitem{d_hajos}
G. A. Dirac, {\it A property of 4-chromatic graphs and some remarks on critical graphs}, J. London Math. Soc. 27 (1952), 85--92.

\bibitem{d_cut}
G. A. Dirac, {\it The structure of $k$-chromatic graphs}, Fund. Math. 40 (1953), 42--55.

\bibitem{dmw}
V, Dujmovi\'{c}, P. Morin, and D. R. Wood, {\it Layered separators in minor-closed graph classes with applications}, J. Combin. Theory Ser. B 127 (2017), 111--147.

\bibitem{dn}
Z. Dvo\v{r}\'{a}k and S. Norin, {\it Islands in minor-closed classes. I. Bounded treewidth and separators}, arXiv:1710.02727.

\bibitem{dy}
Z. Dvo\v{r}\'{a}k and L. Yepremyan, {\it Complete graph immersions and minimum degree}, J. Graph Theory 88 (2018), 211--221.

\bibitem{ekkos}
K. Edwards, D. Y. Kang, J. Kim, S. Oum and P. Seymour, {\it A relative of Hadwiger's conjecture}, SIAM J. Discrete Math. 29 (2015), 2385--2388.

\bibitem{e}
P. Erd\H{o}s, {\it Graph theory and probability}, Canad. J. Math. 11 (1959), 34--38.

\bibitem{ef}
P. Erd\H{o}s and S. Fajtlowicz, {\it On the conjecture of Haj\'{o}s}, Combinatorica 1 (1981), 141--143.

\bibitem{glw}
G. Gauthier, T.-N. Le and P. Wollan, {\it Forcing clique immersions through chromatic number}, European J. Combin. 81 (2019), 98--118.

\bibitem{h}
H. Hadwiger, {\it \"{U}ber eine Klassifikation der Streckenkomplexe}, Vierteljschr. Naturforsch. Ges. Z\"{u}rich 88 (1943), 133--142.

\bibitem{hst}
P. Haxell, T. Szab\'{o} and G. Tardos, {\it Bounded size components--partitions and transversals}, J. Combin. Theory Ser. B 88 (2003), 281--297.

\bibitem{hw}
J. van den Heuvel and D. R. Wood, {\it Improper colourings inspired by Hadwiger's conjecture}, J. London Math. Soc. 98 (2018), 129--148.

\bibitem{km}
K. Kawarabayashi and B. Mohar, {\it A relaxed Hadwiger's conjecture for list colorings}, J. Combin. Theory Ser. B 97 (2007), 647--651.

\bibitem{k1}
A. V. Kostochka, {\it The minimum Hadwiger number for graphs with a given mean degree of vertices}, Metody Diskret. Analiz. 38 (1982), 37--58.

\bibitem{k2}
A. V. Kostochka, {\it Lower bound of the Hadwiger number of graphs by their average degree}, Combinatorica 4 (1984), 307--316.

\bibitem{lm}
F. Lescure and H. Meyniel, {\it On a problem upon configurations contained in graphs with given chromatic number}, Graph Theory in Memory of G.A. Dirac, Ann. Discrete Math. 41 (1989), 325--332.

\bibitem{lmst}
N. Linial, J. Matou\v{s}ek, O. Sheffet and G. Tardos, {\it Graph colouring with no large monochromatic components}, Combin. Probab. Comput. 17 (2008), 577--589.

\bibitem{l_imm_str}
C.-H. Liu, {\it A global decomposition theorem for excluding immersions in graphs with no edge-cut of order three}, arXiv:2006.15694.

\bibitem{lo}
C.-H. Liu and S. Oum, {\it Partitioning $H$-minor free graphs into three subgraphs with no large components}, J. Combin. Theory Ser. B 128 (2018), 114--133.

\bibitem{lw_minor}
C.-H. Liu and D. R. Wood, {\it Clustered coloring of graphs excluding a subgraph and a minor}, arXiv:1905.09495.

\bibitem{lw_layer}
C.-H. Liu and D. R. Wood, {\it Clustered graph coloring and layered treewidth}, arXiv:1905.08969.

\bibitem{lw_topo}
C.-H. Liu and D. R. Wood, {\it Clustered variants of Haj\'{o}s' conjecture}, J. Combin. Theory Ser. B 152 (2022), 27--54.

\bibitem{l_partition}
L. Lov\'{a}sz, {\it On decomposition of graphs}, Stud. Sci. Math. Hung. 1 (1966), 237--238.

\bibitem{ns}
S. Norin, L. Postle and Z.-X Song, {\it Breaking the degeneracy barrier for coloring graphs with no $K_t$ minor}, arXiv:1910.09378.

\bibitem{rsst}
N. Robertson, D. P. Sanders, P. Seymour and R. Thomas, {\it The four-colour theorem}, J. Combin. Theory Ser. B 70 (1997), 2--44.

\bibitem{rst}
N. Robertson, P. Seymour and R. Thomas, {\it Hadwiger's conjecture for $K_6$-free graphs}, Combinatorica 13 (1993), 279--361.

\bibitem{t}
A. Thomason, {\it An extremal function for contractions of graphs}, Math. Proc. Cambridge Philos. Soc. 95 (1984), 261--265.

\bibitem{w_color}
K. Wagner, {\it \"{U}ber eine Eigenschaft der ebenen Komplexe}, Math. Ann. 114 (1937), 570--590.

\bibitem{w_clu}
D. R. Wood, {\it Contractibility and the Hadwiger conjecture}, European J. Combin. 31 (2010), 2102--2109.

\end{thebibliography}
\end{document}